\documentclass[letterpaper, 10 pt, conference]{ieeeconf}  
\IEEEoverridecommandlockouts                              
\usepackage{ieee_packages}

\newtheorem{prop}{Proposition}

\title{\LARGE \bf
Geometric Adaptive Control of Attitude Dynamics on $\SO$\\ with State Inequality Constraints}

\author{Shankar Kulumani, Christopher Poole, and Taeyoung Lee
\thanks{Shankar Kulumani, Christopher Poole, Taeyoung Lee, Mechanical and Aerospace Engineering, George Washington University, Washington DC 20052 {\tt \{skulumani,poolec,tylee\}@gwu.edu}}
\thanks{This research has been supported in part by NSF under the grants CMMI-1243000, CMMI-1335008, and CNS-1337722.}
}

\begin{document}
\allowdisplaybreaks

\maketitle
\thispagestyle{empty}
\pagestyle{empty}

\begin{abstract}
This paper presents a new geometric adaptive control system with state inequality constraints for the attitude dynamics of a rigid body. 
The control system is designed such that the desired attitude is asymptotically stabilized, while the controlled attitude trajectory avoids undesired regions defined by an inequality constraint. 
In addition, we develop an adaptive update law that enables attitude stabilization in the presence of unknown disturbances. 
The attitude dynamics and the proposed control systems are developed on the special orthogonal group such that singularities and ambiguities of other attitude parameterizations, such as Euler angles and quaternions are completely avoided. 
The effectiveness of the proposed control system is demonstrated through numerical simulations and experimental results.
\end{abstract}

\section{Introduction}\label{sec:intro}

Rigid body attitude control is an important problem for aerospace vehicles, ground and underwater vehicles, as well as robotic systems~\cite{hughes2004,wertz1978}.
One distinctive feature of the attitude dynamics of rigid bodies is that it evolves on a nonlinear manifold.
The three-dimensional special orthogonal group, or \( \SO \), is the set of \( 3 \times 3 \) orthogonal matrices whose determinant is one.
This configuration space is non-Euclidean and yields unique stability properties which are not observable on a linear space.
For example, it is impossible to achieve global attitude stabilization using continuous time-invariant feedback~\cite{bhat2000}.



Attitude control is typically studied using a variety of attitude parameterizations, such as Euler angles or quaternions~\cite{shuster1993}.
All attitude parameterizations fail to represent the nonlinear configuration space both globally and uniquely~\cite{chaturvedi2011a}.
For example, minimal attitude representations, such as Euler angle sequences or modified Rodriguez parameters, suffer from singularities.
These attitude representations are not suitable for large angular slews.
Quaternions do not have singularities but they double cover the special orthogonal group.
As a result, any physical attitude is represented by a pair of antipodal quaternions on the three-sphere.
During implementation, the designer must carefully resolve this non-unique representation in quaternion based attitude control systems to avoid undesirable unwinding behavior~\cite{bhat2000}.

Many physical rigid body systems must perform large angular slews in the presence of state constraints.
For example, autonomous spacecraft or aerial systems are typically equipped with sensitive optical payloads, such as infrared or interferometric sensors.
These systems require retargeting while avoiding direct exposure to sunlight or other bright objects.
The removal of constrained regions from the rotational configuration space results in a \textit{nonconvex} region.
The attitude control problem in the absence of constraints has been extensively studied~\cite{bullo2004,MayTeePaCC11,LEEITAC15}.
However, the attitude control problem in the presence of constraints has received much less attention.

Several approaches have been developed to treat the attitude control problem in the presence of constraints.
A conceptually straightforward approach is used in~\cite{hablani1999} to determine feasible attitude trajectories prior to implementation.
The algorithm determines an intermediate point such that an unconstrained maneuver can be calculated for each subsegment.
Typically, an optimal or easily implementable on-board control scheme for attitude maneuvers is applied to maneuver the vehicle along these segments.
In this manner it is possible to accomplish constraint avoidance by linking several intermediary unconstrained maneuvers.
While this method is conceptually simple, it is difficult to generalize for an arbitrary number of constraints.
In addition, this approach is only applicable to problems where the selection of intermediate points are computationally feasible.

The approach in~\cite{frazzoli2001} involves the use of randomized motion planning algorithms to solve the constrained attitude control problem.
A graph is generated consisting of vertices from an initial attitude to a desired attitude. 
A random iterative search is conducted to determine a path through a directed graph such that a given cost functional is minimized.
The random search approach can only stochastically guarantee attitude convergence as it can be shown that as the number of vertices in the graph grow, the probability of nonconvergence goes to zero.
However, the computational demand grows as the size of the graph is increased. 
As a result, random search approaches are ill-suited to on-board implementation or in scenarios that require agile maneuvers.

Model predictive control for spacecraft attitude dynamics is studied in~\cite{guiggiani2014,kalabic2014,gupta2015}.
These methods rely on linear or non-linear state dynamics to repeatedly solve a finite-time optimal control problem.
As a result, model predictive control methods are also computational expensive and apply direct optimization methods to solve the necessary conditions for optimality.
Therefore these methods are complicated to implement and not applicable for real-time control applications.
  
Artificial potential functions are commonly used to handle kinematic constraints for a wide range of problems in robotics~\cite{rimon1992}.
The goal is the design of attractive and repulsive terms which drive the system toward or away from a certain state, respectively.
The superposition of the these functions allows one to apply standard feedback control schemes for stabilization and tracking.
More specifically, artificial potential functions have previously been applied to the spacecraft attitude control problem in~\cite{lee2011b,mcinnes1994}.
However, both of these approaches were developed using attitude parameterizations, namely Euler angles and quaternions, and as such, they are limited by the singularities of minimal representations or the ambiguity of quaternions.

This paper is focused on developing an adaptive attitude control scheme in the presence of attitude inequality constraints on \(\SO\).
We apply a potential function based approach developed directly on the nonlinear manifold \(\SO\). 
By characterizing the attitude both globally and uniquely on \(\SO\), our approach avoids the issues of attitude parameterizations, such as kinematic singularities and ambiguities, and is geometrically exact. 
A configuration error function on \(\SO\) with a logarithmic barrier function is proposed to avoid constrained regions. 
Instead of calculating a priori trajectories, as in the geometric and randomized approaches, our approach results in a closed-loop attitude control system. 
This makes it ideal for on-board implementation on UAV or spacecraft systems. 
In addition, unlike previous approaches our control system can handle an arbitrary number of constrained regions without modification.

Furthermore, we formulate an adaptive update law to enable attitude convergence in the presence of uncertain disturbances. 
The stability of the proposed control systems is verified via mathematically rigorous Lyapunov analysis on $\SO$.  
In short, the proposed attitude control system in the presence of inequality constraints is computationally efficient and able to handle uncertain disturbances. 
The effectiveness of this approach is illustrated via numerical simulation and experimental results.

\section{Problem Formulation}\label{sec:prob_form}
\subsection{Attitude Dynamics}\label{sec:att_dyn}
Consider the attitude dynamics of a rigid body. 
We define an inertial reference frame and a body frame whose origin is at the center of mass and aligned with the principle directions of the body. 
The configuration manifold of the attitude dynamics is the special orthogonal group:
\begin{align*}
	\SO = \{R\in\R^{3\times 3}\,|\, R^TR=I,\;\mathrm{det}[R]=1\} ,
\end{align*}
where a rotation matrix $R\in\SO$ represents the transformation of the representation of a vector from the body-fixed frame to the inertial reference frame. 
The equations of motion are given by
\begin{gather}
	J\dot\Omega + \Omega\times J\Omega = u+W(R,\Omega)\Delta ,\label{eqn:Wdot}\\
	\dot R = R\hat\Omega ,\label{eqn:Rdot}
\end{gather}
where $J\in\R^{3\times 3}$ is the inertia matrix, and $\Omega\in\R^{3}$ is the angular velocity represented with respect to the body-fixed frame. 
The control moment is denoted by $u\in\R^{3}$, and it is expressed with respect to the body-fixed frame. 
We assume that the external disturbance is expressed by $W(R,\Omega)\Delta$, where $W(R,\Omega):\SO\times\R^{3}\rightarrow \R^{3\times p}$ is a known function of the attitude and the angular velocity.
The disturbance is represented by $\Delta\in\R^{p}$ and is an unknown, but fixed uncertain parameter.
In addition, we assume that a bound on \( W(R, \Omega) \text{ and } \Delta \) is known and given by
\begin{equation}
	\norm{W} \leq B_W , \quad \norm{\Delta} \leq B_\Delta \,. \label{eqn:W_bound}
\end{equation}
This form of uncertainty enters the system dynamics through the input channel and as a result is referred to as a matched uncertainty. 
While this form of uncertainty is easier than the unmatched variety many physically realizable disturbances may be modeled in this manner.
For example, orbital spacecraft are subject to gravity gradient torques caused by the non-spherical distribution of mass of both the spacecraft and central gravitational body.
This form of disturbance may be represented as a body fixed torque on the vehicle.
In addition, for typical scenarios, where the spacecraft is significantly smaller than the orbital radius, the disturbance torque may be assumed constant over short time intervals.

In~\refeqn{Rdot}, the \textit{hat} map $\wedge :\R^{3}\rightarrow\so$ represents the transformation of a vector in $\R^{3}$ to a $3\times 3$ skew-symmetric matrix such that $\hat x y = x\times y$ for any $x,y\in\R^{3}$~\cite{bullo2004}. 
More explicitly, 
\begin{align*}
\hat x = \begin{bmatrix} 0 & -x_3 & x_2 \\ x_3 & 0 & -x_1 \\ -x_2 & x_1 & 0\end{bmatrix},
\end{align*}
for $x=[x_1,x_2,x_3]^T\in\R^{3}$. 
The inverse of the hat map is denoted by the \textit{vee} map $\vee:\so\rightarrow\R^{3}$. 
Several properties of the hat map are summarized as
\begin{gather}
    x\cdot \hat y z = y\cdot \hat z x,\quad \hat x\hat y z = (x\cdot z) y - (x\cdot y ) z\label{eqn:STP},\\
    \widehat{x\times y} = \hat x \hat y -\hat y \hat x = yx^T-xy^T,\label{eqn:hatxy}\\
    \tr{A\hat x }=\frac{1}{2}\tr{\hat x (A-A^T)}=-x^T (A-A^T)^\vee,\label{eqn:hat1}\\
    \hat x  A+A^T\hat x=(\braces{\tr{A}I_{3\times 3}-A}x)^{\wedge},\label{eqn:xAAx}\\
R\hat x R^T = (Rx)^\wedge,\quad 
R(x\times y) = Rx\times Ry\label{eqn:RxR}
\end{gather}
for any $x,y,z\in\R^{3}$, $A\in\R^{3\times 3}$ and $R\in\SO$. 
Throughout this paper, the dot product of two vectors is denoted by $x\cdot y = x^T y$ for any $x,y\in\R^n$ and the maximum eigenvalue and the minimum eigenvalue of $J$ are denoted by $\lambda_M$ and $\lambda_m$, respectively. 
The 2-norm of a matrix \( A \) is denoted by \( \norm{A} \), and its Frobenius norm is denoted by \( \norm{A} \leq \norm{A}_F = \sqrt{\tr{A^T A}} \leq \sqrt{\text{rank}(A)} \norm{A} \).

\subsection{State Inequality Constraint}

The two-sphere is the manifold of unit-vectors in \( \R^3 \) such that \( \Sph^2 = \{ q \in \R^3 \,  \vert \, \norm{q} = 1 \}\).
We define \( r \in \Sph^2 \) to be a unit vector from the mass center of the rigid body along a certain direction and it is represented with respect to the body-fixed frame.
For example, \( r \) may represent the pointing direction of an on-board optical sensor.
We define \( v \in \Sph^2 \) to be a unit vector from the mass center of the rigid body toward an undesired pointing direction and represented in the inertial reference frame.
For example, \( v \) may represent the inertial direction of a bright celestial object or the incoming direction of particles or other debris.
It is further assumed that optical sensor has a strict non-exposure constraint with respect to the celestial object.
We formulate this hard constraint as
\begin{align}
	r^T R^T v \leq \cos \theta , \label{eqn:constraint}
\end{align}
where we assume \( \ang{0} \leq \theta \leq \ang{90}  \) is the required minimum angular separation between \( r \) and \( R^T v \). 

The objective is to a determine a control input \( u \) that stabilizes the system from an initial attitude \( R_0 \) to a desired attitude \( R_d \) while ensuring that~\cref{eqn:constraint} is always satisfied.

\section{Attitude Control on $\SO$ with Inequality Constraints}

The first step in designing a control system on a nonlinear manifold \( \Q \) is the selection of a proper configuration error function. 
This configuration error function, \( \Psi : \Q \times \Q \to \R \), is a smooth and proper positive definite function that measures the error between the current configuration and a desired configuration. 
Once an appropriate configuration error function is chosen, one can then define a configuration error vector and a velocity error vector in the tangent space \( \mathsf{T}_q \Q \) through the derivatives of \( \Psi \)~\cite{bullo2004}. 
With the configuration error function and vectors the remaining procedure is analogous to nonlinear control design on Euclidean vector spaces. 
One chooses control inputs as functions of the state through a Lyapunov analysis on \Q.

To handle the attitude inequality constraint, we propose a new attitude configuration error function. 
More explicitly, we extend the trace form used in~\cite{bullo2004,LeeITCST13} for attitude control on \(\SO\) with the addition of a logarithmic barrier function. 
Based on the proposed configuration error function,  nonlinear geometric attitude controllers are constructed. 
A smooth control system is first developed assuming that there is  no disturbance, and then it is extended to include an adaptive update law for stabilization in the presence of unknown disturbances. 
The proposed attitude configuration error function and several properties are summarized as follows.

\begin{prop}[Attitude Error Function] \label{prop:config_error}
Define an attitude error function \( \Psi : \SO \to \R \), an attitude error vector \( e_R \in \R^3 \), and an angular velocity error vector \( e_\Omega \in \R^3 \) as follows:
\begin{gather}
	\Psi(R) = A(R) B(R) , \label{eqn:psi} \\
	e_R = e_{R_A} B(R) + A(R) e_{R_B} , \label{eqn:eR} \\
	e_\Omega = \Omega , \label{eqn:eW}
\end{gather}
with
\begin{gather}
	A(R) = \frac{1}{2} \tr{G \left( I - R_d^T R\right)} , \label{eqn:A} \\
	B(R) = 1 - \frac{1}{\alpha} \ln \left( \frac{\cos \theta -  r^T R^T v}{1 + \cos \theta}\right) . \label{eqn:B} \\
	e_{R_A} = \frac{1}{2} \parenth{G R_d^T R - R^T R_d G} ^ \vee , \label{eqn:eRA} \\
	e_{R_B} = \frac{ \left( R^T v\right)^\vee r}{\alpha \left(r^T R^T v - \cos \theta \right)} . \label{eqn:eRB} 
\end{gather}	
where \( \alpha \in \R \) is defined as a positive constant and the matrix \( G \in \R^{3 \times 3} \) is defined as a diagonal matrix matrix for distinct, positive constants \( g_1, g_2, g_3 \in \R \).
Then, the following properties hold
\begin{enumerate}[(i)]
	\item \label{item:prop_psi_psd} \(\Psi\) is positive definite about \( R = R_d\)
	\item \label{item:prop_era}The variation of \( A(R) \) with respect to a variation of \( \delta R = R \hat{\eta} \) for \( \eta \in \R^3 \) is given by
	\begin{align}
		\dirDiff{A}{R} &= \eta \cdot e_{R_A} .
	\end{align}
	\item \label{item:prop_erb} The variation of \( B(R) \) with respect to a variation of \( \delta R = R \hat{\eta} \) for \( \eta \in \R^3 \) is given by
	\begin{align}
		\dirDiff{B}{R} &= \eta \cdot e_{R_{B}} .
	\end{align}
	\item \label{item:prop_crit}The critical points of \( \Psi \) are $R_d$, and $R_d \exp(\pi \hat{s})$ for $s \in \braces{e_1, e_2, e_3}$ satisfying $R^T v = \pm r$.
	\item \label{item:prop_era_upbound}An upper bound of \( \norm{e_{R_A}} \) is given as:
	\begin{align}
		\norm{e_{R_A}}^2 \leq \frac{A(R)}{b_1} , \label{eqn:psi_lower_bound}
	\end{align}
	where the constant \( b_1 \) is given by \( b_1 = \frac{h_1}{h_2 + h_3} \) for 
	\begin{align*}
		h_1 &= \min\braces{g_1 + g_2, g_2 + g_3 , g_3 + g_1} ,\\
		h_2 &= \min\braces{\parenth{g_1 -g_2}^2,\parenth{g_2 -g_3}^2 , \parenth{g_3 -g_1}^2} ,\\
		h_3 &= \min\braces{\parenth{g_1 + g_2}^2, \parenth{g_2 + g_3}^2 , \parenth{g_3 + g_1}^2}.		
	\end{align*}
\end{enumerate}
\end{prop}
\begin{proof}
See~\Cref{proof:config_error}.
\end{proof}

\begin{figure} 
	\centering 
	\begin{subfigure}[htbp]{0.45\columnwidth} 
		\includegraphics[width=\columnwidth]{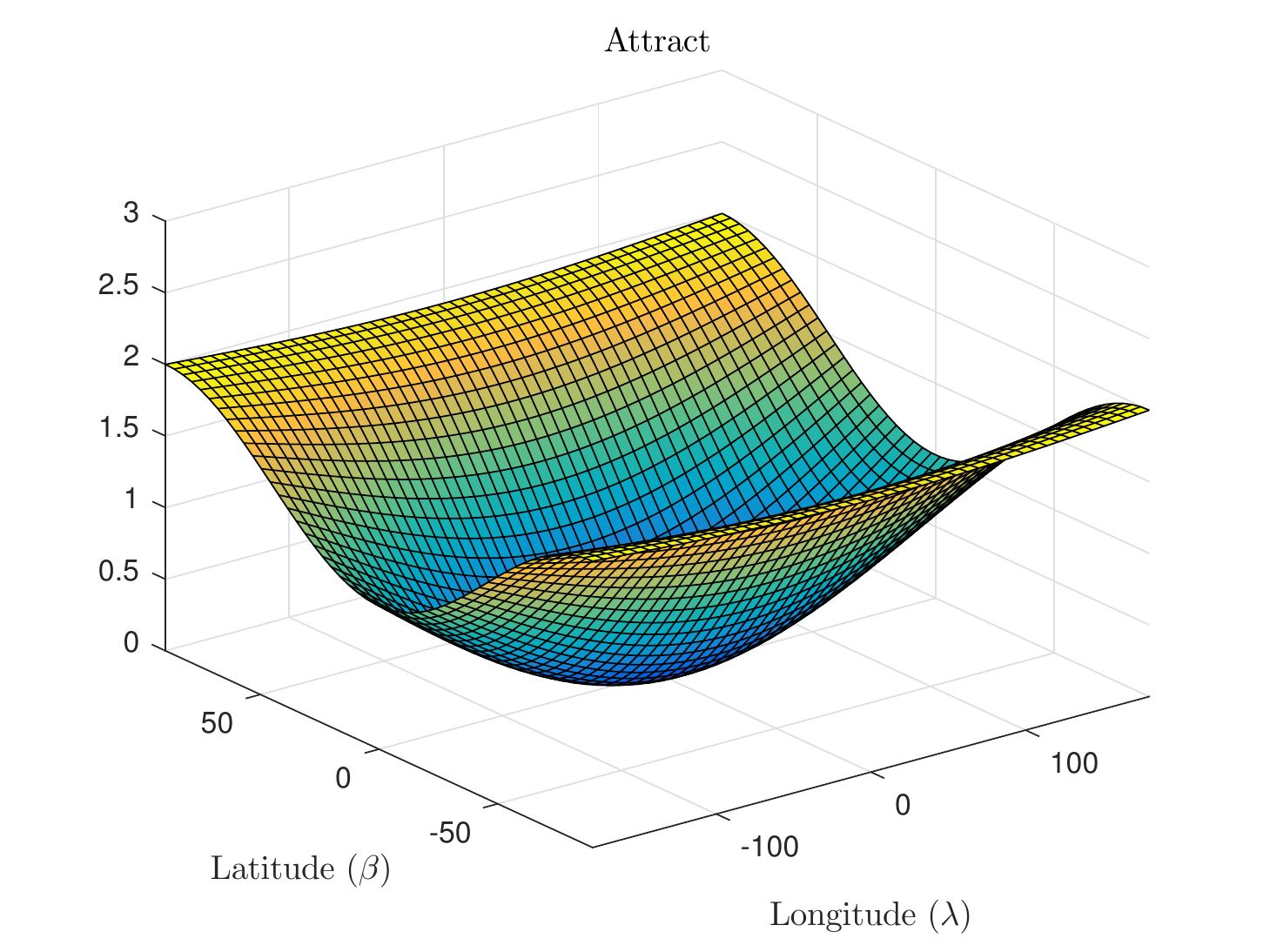} 
		\caption{Attractive \( A(R) \) } \label{fig:attract_error} 
	\end{subfigure}~ 
	\begin{subfigure}[htbp]{0.45\columnwidth} 
		\includegraphics[width=\columnwidth]{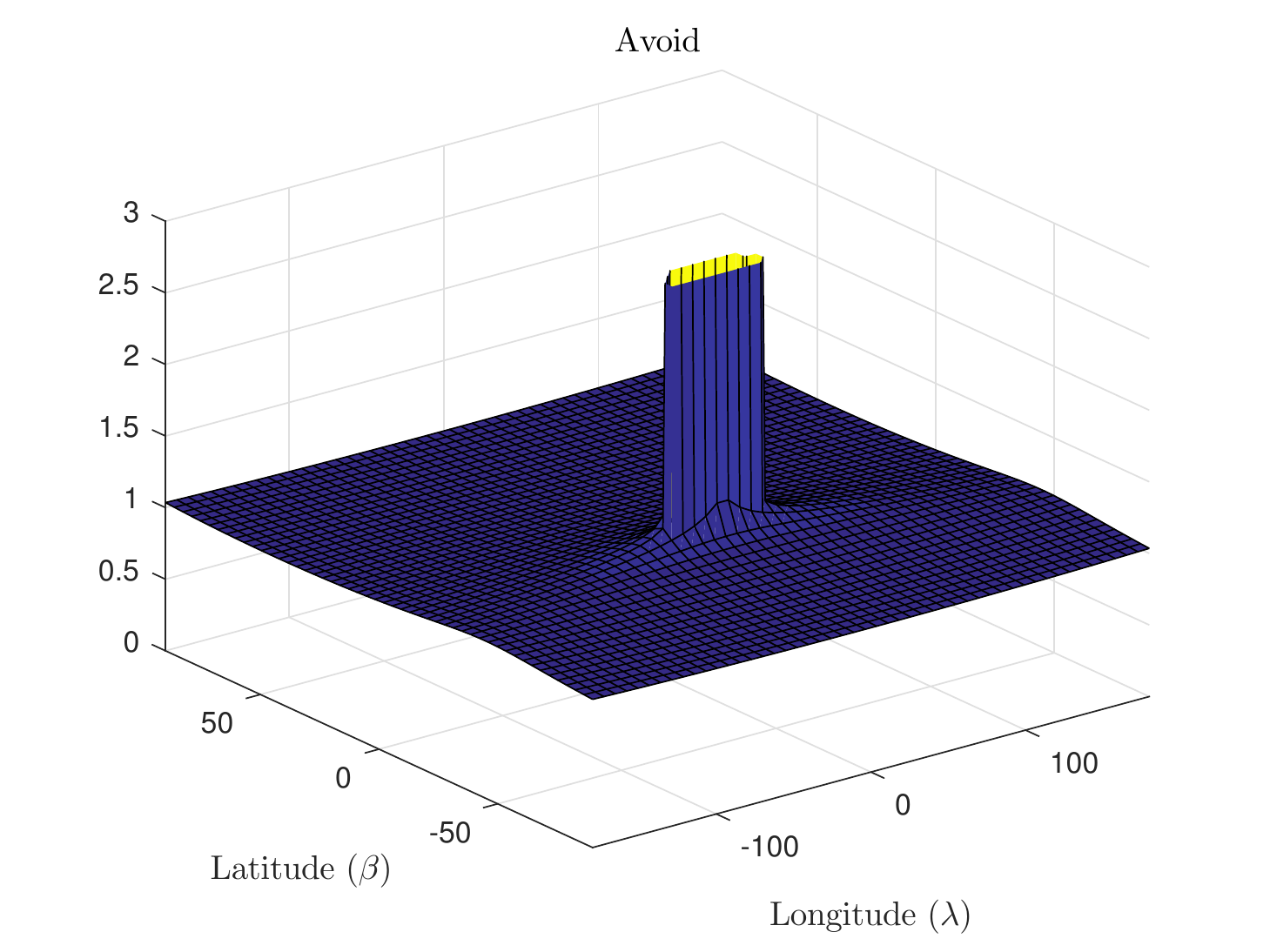} 
		\caption{Repulsive \( B(R) \)} \label{fig:avoid_error} 
	\end{subfigure}
	\centering
	\begin{subfigure}[htbp]{0.45\columnwidth} 
		\includegraphics[width=\columnwidth]{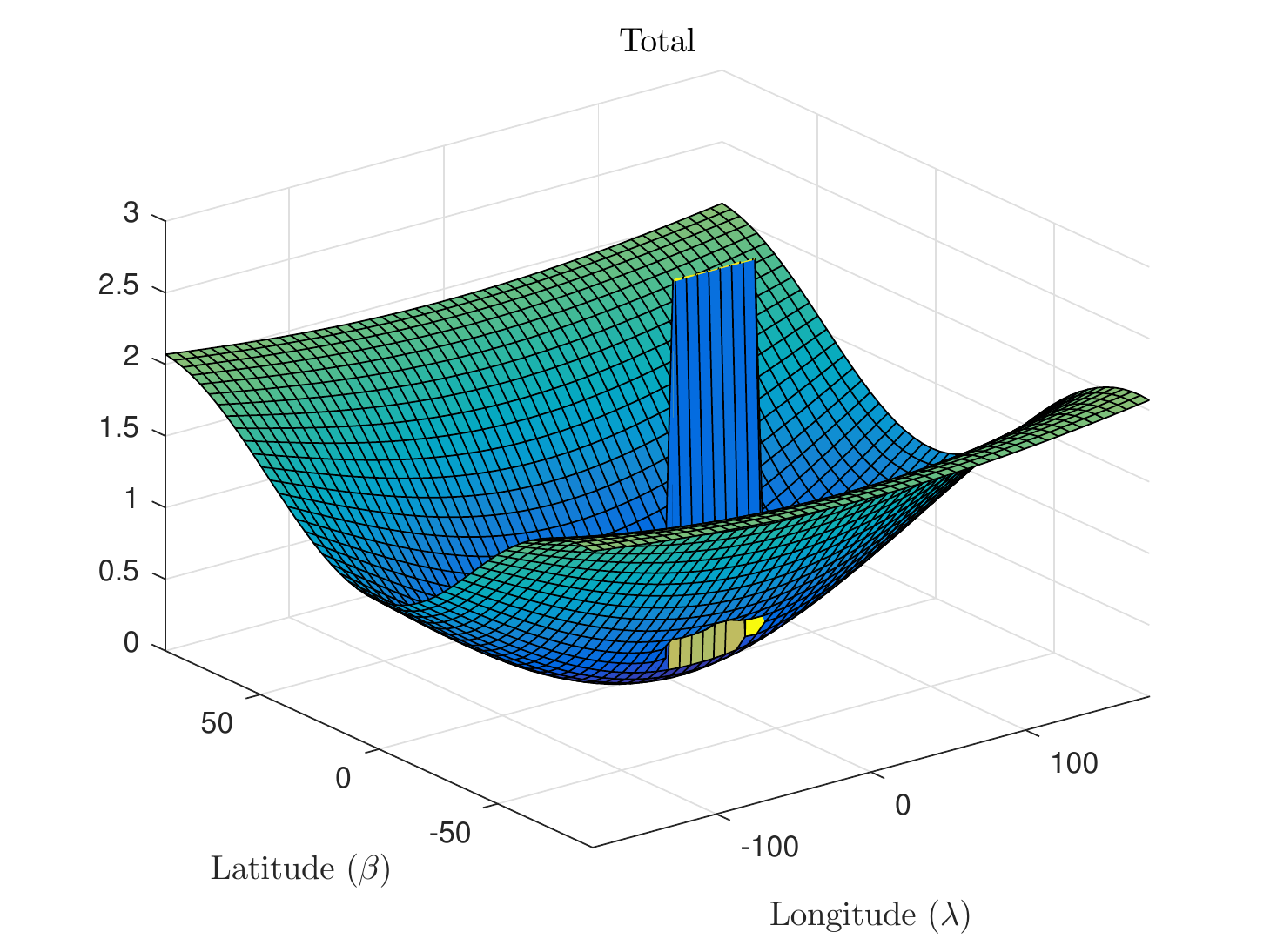} 
		\caption{Configuration \( \Psi \)} \label{fig:combined_error} 
	\end{subfigure}
	\caption{Configuration error function visualization}
	\label{fig:config_error} 
\end{figure}

\Cref{eqn:psi} is composed of an attractive term, \( A (R) \) toward the desired attitude, and a repulsive term, \( B(R) \) away from the undesired direction \( R^T v \).
In order to visualize the attitude error function on \( \SO \) we utilize a spherical coordinate representation.
Recall that the spherical coordinate system represents the position of a point relative to an origin in terms of a radial distance, azimuth, and elevation.
This coordinate system is commonly used to define locations on the Earth in terms of a latitude and longitude.
Similarly, the positions of celestial objects are defined on the celestial sphere in terms of right ascension and declination. 
We apply this concept and parametrize the rotation matrix \( R \in \SO \) in terms of the spherical angles \( \SI{-180}{\degree} \leq \lambda \leq \SI{180}{\degree}  \) and \( \SI{-90}{\degree} \leq \beta \leq \SI{90}{\degree} \). 
Using the elementary Euler rotations the rotation matrix is now defined as \( R = \exp( \lambda \hat{e}_2) \exp( \beta \hat{e}_3) \).
We iterate over the domains of \( \lambda\) and \(\beta\) in order to rotate the body-fixed vector \( r \) throughout the two-sphere \( \S^2 \).
Applying this method,~\cref{fig:config_error} allows us to visualize the error function on \( \SO \).
The attractive error function, given by~\cref{eqn:A}, has been previously used for attitude control on \(\SO\).
The potential well of \( A(R)\) is illustrated in~\cref{fig:attract_error}, where the desired attitude lies at the minimum of \( A(R) \).

To incorporate the state inequality constraints we apply a logarithmic barrier term.
Barrier functions are typically used in optimal control and motion planning applications.
A visualization of the configuration error function is presented in~\cref{fig:avoid_error} which shows that as the boundary of the constraint is neared, or \( r^T R^T v \to \cos \theta \), the barrier term increases, \( B \to \infty\).
We use the scale factor~\(\frac{1}{1+\cos \theta} \) to ensure that \( \Psi \) remains positive definite.
The logarithmic function is popular as it quickly decays away from the constraint boundary.
The positive constant \( \alpha \) serves to shape the barrier function.
As \( \alpha \) is increased the impact of \( B(R) \) is reduced away from the constraint boundary. 
The superposition of the attractive and repulsive functions is shown in~\cref{fig:combined_error}.
The control system is defined such that the attitude trajectory follows the negative gradient of \( \Psi \) toward the minimum at \( R = R_d \), while avoiding the constrained region.

While~\cref{eqn:B} represents a single inequality constraint given as~\cref{eqn:constraint}, it is readily generalized to multiple constraints of an arbitrary form. 
For example, the configuration error function can be formulated as $\Psi=A[1+\sum_i C_i]$, where $C_i$ has the form of $C_i=B-1$ for the $i$-th constraint. 
In this manner, one may enforce multiple state inequality constraints, and we later demonstrate this through numerical simulation. 

\begin{prop}[Error Dynamics]\label{prop:error_dyn}
	The attitude error dynamics for \( \Psi, e_R, e_\Omega \) satisfy 
	\begin{gather}
		\diff{}{t} \parenth{\Psi} = e_R \cdot e_\Omega , \label{eqn:psi_dot}\\
		\diff{}{t} \parenth{e_R} = \dot{e}_{R_A} B(R) + e_{R_A} \dot{B}(R) + \dot{A}(R)e_{R_B} + A \dot{e}_{R_B} , \label{eqn:eR_dot}\\
		\diff{}{t} \parenth{e_{R_A}} = E(R, R_d) e_\Omega , \label{eqn:eRA_dot} \\
		\diff{}{t} \parenth{e_{R_B}} = F(R) e_\Omega , \label{eqn:eRB_dot} \\
		\diff{}{t} \parenth{A(R)} = e_{R_A} \cdot e_\Omega , \label{eqn:A_dot} \\
		\diff{}{t} \parenth{B(R)} = e_{R_B} \cdot e_\Omega , \label{eqn:B_dot} \\
		\diff{}{t} \parenth{e_\Omega} = J^{-1} \parenth{-\Omega \times J \Omega + u + W(R, \Omega) \Delta} , \label{eqn:eW_dot}
	\end{gather}
	where the matrices \(E(R,R_d), F(R) \in \R^{3\times3} \) are given by
	\begin{gather}
		E(R,R_d) = \frac{1}{2} \parenth{\tr{R^T R_d G}I - R^T R_d G} , \label{eqn:E} \\
		F(R) = \frac{1}{\alpha \parenth{r^T R^T v - \cos \theta}} \left[\parenth{v^T R r} I - R^T v r^T + \right. \nonumber \\
		\left. \frac{R^T \hat{v} R r v^T R \hat{r}}{\parenth{r^T R^T v - \cos \theta}}\right] . \label{eqn:F}
	\end{gather}
\end{prop}
\begin{proof}
See~\Cref{proof:error_dyn}.
\end{proof}

\subsection{Attitude Control without Disturbance}
We introduce a nonlinear geometric controller for the attitude stabilization of a rigid body.
We first assume that there is no disturbance, i.e., \( \Delta = 0 \).
\begin{prop}[Attitude Control]\label{prop:att_control}
	Given a desired attitude command \( \parenth{R_d, \Omega_d = 0} \), which satisfies the constraint~\cref{eqn:constraint}, and positive constants \( k_R, k_\Omega \in \R \) we define a control input \( u \in \R^3 \) as follows
	\begin{gather}
		u = -k_R e_R - k_\Omega e_\Omega + \Omega \times J \Omega . \label{eqn:nodist_control}
	\end{gather}
	Then the zero equilibrium of the attitude error is asymptotically stable, and the inequality constraint is satisfied.
\end{prop}

\begin{proof}
See~\Cref{proof:att_control}.
\end{proof}

This proposition only guarantees that the attitude error vector \( e_R \) asymptotically converges to zero.
However, this does not necessarily imply that \( R \to R_d \) as \( t \to \infty \), since there are at most three additional critical points of \( \Psi \) where \( e_R = 0 \) and \( R^T v = \pm r\).
At an undesired equilibrium \( R = \exp{ \parenth{\pi \hat{e}_i}} R_d \) and \( e_\Omega =0 \).
However, we can show that these undesired equilibrium points are unstable in the sense of Lyapunov~\cite{LeeITCST13}.
As a result, we can claim that the desired equilibrium \( R = R_d \text{ and } e_\Omega = 0 \) is almost globally asymptotically stable, which means that the set of initial conditions that do not converge to the desired attitude has zero Lebesgue measure.
 
\subsection{Adaptive Control}
We extend the results of the previous section with the addition of a fixed but unknown disturbance \( \Delta \).
This scenario is typical of many mechanical systems and represents unmodeled dynamics or external moments acting on the system.
For example, Earth orbiting spacecraft typically experience a torque due to a gravitational gradient.
Aerial vehicles will similarly experience external torques due to air currents or turbulence.
An adaptive control system is introduced to asymptotically stabilize the system to a desired attitude while ensuring that state constraints are satisfied. 

\begin{prop}[Bound on \( \dot{e}_R \)]\label{prop:eR_dot_bound}
Consider a domain \( D \) about the desired attitude defined as
\begin{align}
	D = \braces{R \in \SO \vert \Psi < \psi < h_1, r^T R^T v < \beta < \cos \theta}. \label{eqn:domain}
\end{align}
Then the following statements hold:
\begin{enumerate}[(i)]
	\item \label{item:prop_eR_dot_bound_AB} Upper bounds of \( A(R) \) and \( B(R) \) are given by
	\begin{gather}
		\norm{A} < c_A  , \quad \norm{B} < c_B . \label{eqn:AB_bound}
	\end{gather}
	\item \label{item:prop_eR_dot_bound_EF} Upper bounds of \( E(R,R_d) \) and \( F(R) \) are given by
	\begin{gather}
		\norm{E} \leq \frac{1}{\sqrt{2}} \tr{G}  , \label{eqn:E_bound} \\
		\norm{F} \leq \frac{\parenth{\beta^2 + 1}\parenth{\beta - \cos \theta}^2 + 1 + \beta^2 \parenth{\beta^2-2}}{\alpha^2 \parenth{\beta-\cos \theta}^4} . \label{eqn:F_bound}
	\end{gather}
	\item Upper bounds of the attitude error vectors \( e_{R_A} \) and \( e_{R_B} \) are given by
	\begin{gather}
		\norm{e_{R_A}} \leq \sqrt{\frac{\psi}{b_1}}, \label{eqn:eRA_bound} \\
		\norm{e_{R_B}} \leq \frac{\sin\theta}{\alpha \parenth{\cos \theta - \beta}}. \label{eqn:eRB_bound}
	\end{gather}
\end{enumerate}
These results are combined to yield a maximum upper bound of the time derivative of the attitude error vector \( \dot{e}_R \) as
\begin{gather}
	\norm{\dot{e}_R} \leq H \norm{e_\Omega} ,\label{eqn:eR_bound}
\end{gather}
where  \( H \in \R \) is defined as
\begin{gather}
	H = \norm{B} \norm{E} + 2 \norm{e_{R_A}} \norm{e_{R_B}} + \norm{A}\norm{F}. \label{eqn:H}
\end{gather}
\end{prop}

\begin{proof}
See~\Cref{proof:eR_dot_bound}.
\end{proof}
\begin{prop}[Adaptive Attitude Control]\label{prop:adaptive_control}
Given  a desired attitude command \( (R_d, \Omega_d = 0 )\) and positive constants \( k_R, k_\Omega, k_\Delta, c \in \R \), we define a control input \( u \in \R^3\) and an adaptive update law for the estimated uncertainty \( \bar{\Delta} \) as follows:
\begin{align}
	u &= - k_R e_R - k_\Omega e_\Omega + \Omega \times J \Omega - W \bar{\Delta} , \label{eqn:adaptive_control} \\
	\dot{\bar{\Delta}} &= k_\Delta W^T \parenth{e_\Omega + c e_R} . \label{eqn:delta_dot}
\end{align}
If \( c \) is chosen such that
\begin{gather}
	0 < c < \frac{4 k_R k_\Omega}{k_\Omega^2 + 4 k-R \lambda_M H} , \label{eqn:c_bound}
\end{gather}
  the zero equilibrium of the error vectors is stable in the sense of Lyapunov. Furthermore, $e_R,e_\Omega\rightarrow 0$ as $t\rightarrow\infty$, and $\bar\Delta$ is uniformly bounded.
\end{prop}
\begin{proof}
See~\Cref{proof:adaptive_control}.
\end{proof}

Nonlinear adaptive controllers have been developed for attitude stabilization in terms of modified Rodriguez parameters and quaternions, as well as attitude tracking in terms of Euler angles. 
The proposed control system is developed on \(\SO\) and avoids the singularities of Euler angles and Rodriguez parameters while incorporating state inequality constraints. 
In addition, the unwinding and double coverage ambiguity of quaternions are also completely avoided. 
The control system handles uncertain disturbances while avoiding constrained regions.

Compared to the previous work on constrained attitude control, we present a geometrically exact control system without parameterizations.
In addition, we incorporate state inequality constraints for the first time on \( \SO \).
The presented control system is computed in real-time and offers significant computational advantages over previous iterative methods. 
In addition, the riguous mathematical proof guarantees stability.
\section{Numerical Examples}
We demonstrate the performance of the proposed control system via numerical simulation.
The inertia tensor of a rigid body is given as
\begin{gather*}
	J = \begin{bmatrix}
	\num{5.57e-3} & \num{6.17e-5} & \num{-2.50e-5} \\
	\num{6.17e-5} & \num{5.57e-3} & \num{1.00e-5} \\
	\num{-2.50e-5} & \num{1.00e-5} & \num{1.05e-2}
	\end{bmatrix} \si{\kilo\gram\meter\squared} .
\end{gather*} 
The control system parameters are chosen as
\begin{gather*}
	G = \text{diag} [0.9,1.1,1.0], \quad k_R = 0.4 , \quad	k_\Omega = 0.296 ,\\
	c = 1.0 , \quad k_\Delta = 0.5 , \quad \alpha = 15 .
\end{gather*}
A body fixed sensor is defined as \(r = [1,0,0]\), while multiple inequality constraints are defined in~\Cref{tab:constraints}.
The simulation parameters are chosen to be similar to those found in~\cite{lee2011b}, however we increase the size of the constraint regions to create a more challenging scenario for the control system.

The initial state is defined as \(R_0 =  \exp(\ang{225} \times \frac{\pi}{180} \hat{e}_3), \Omega_0 = 0\). 
The desired state is \( R_d = I,\Omega_d = 0\).
\begin{table}
\caption{Constraint Parameters~\label{tab:constraints}}
\begin{center}\begin{tabular}{lc}
Constraint Vector (\( v \)) & Angle (\( \theta \)) \\ \hline \hline 
\([0.174,\,-0.934,\, -0.034]^T\) & \ang{40} \\ \hline 
\([0 ,\, 0.7071 ,\, 0.7071]^T\) & \ang{40} \\ \hline 
\([-0.853 ,\, 0.436 ,\, -0.286]^T\) & \ang{40} \\ \hline 
\([-0.122 ,\,-0.140,\, -0.983]^T\) & \ang{20}\end{tabular} 
\end{center}
\end{table}
We show simulation results for the system stabilizing about the desired attitude with and without the adaptive update law from~\Cref{prop:adaptive_control}.
We assume a fixed disturbance of \(\Delta = \begin{bmatrix} 0.2 & 0.2 & 0.2 \end{bmatrix}^T \), with the function \( W(R,\Omega) = I \).
This form is equivalent to an integral control term which penalizes deviations from the desired configuration.
The first term of~\cref{eqn:delta_dot} has the effect of increasing the proportional gain of the control system, since the time derivative of the attitude error vector, \( \dot{e}_{R} \), is linear with respect to the angular velocity error vector \( e_\Omega\).
\begin{figure} 
	\centering 
	\begin{subfigure}[htbp]{0.5\columnwidth} 
		\includegraphics[width=\columnwidth]{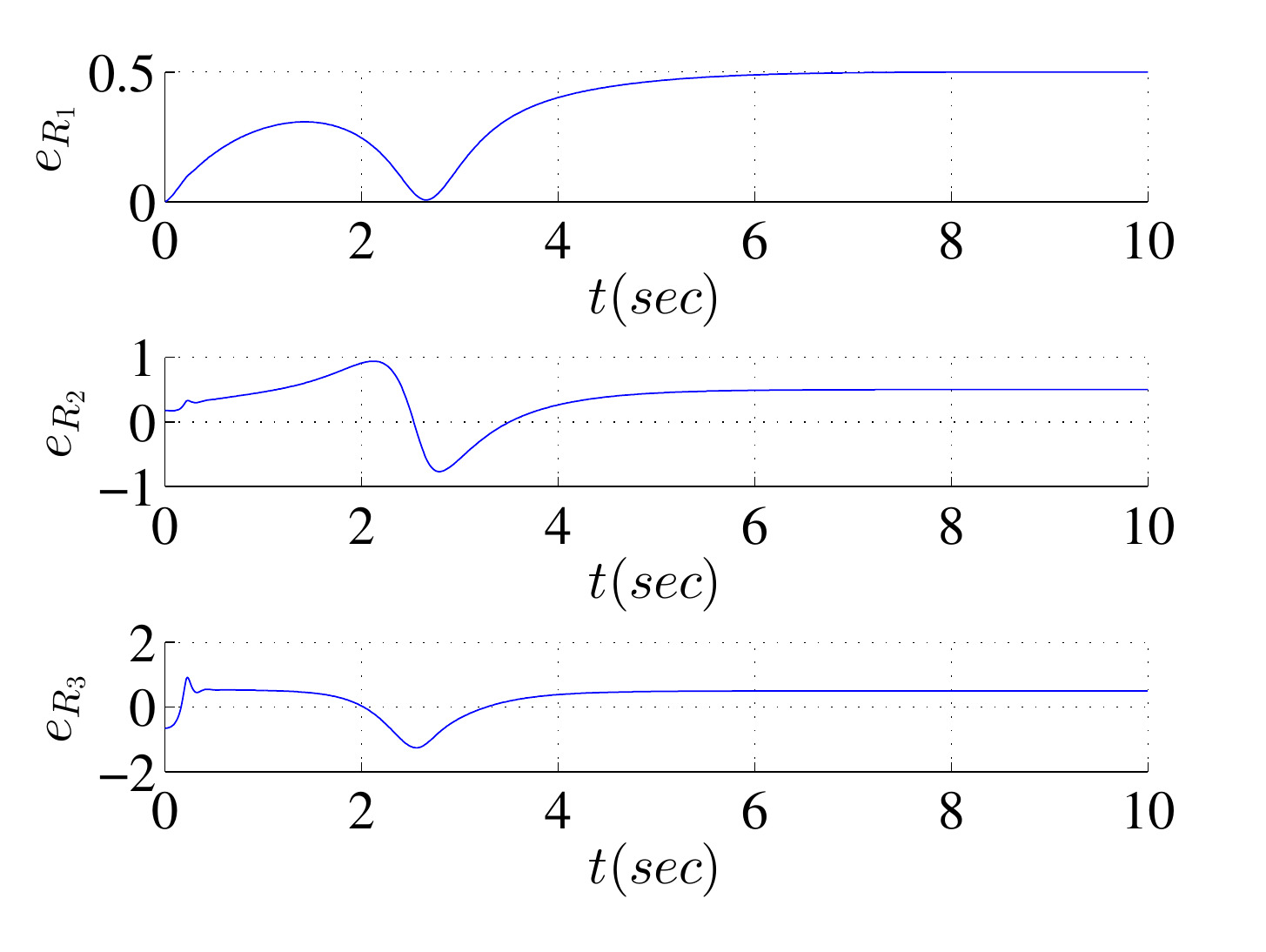} 
		\caption{Attitude error vector \(e_R\) } \label{fig:eR_con} 
	\end{subfigure}~ 
	\begin{subfigure}[htbp]{0.5\columnwidth} 
		\includegraphics[width=\columnwidth]{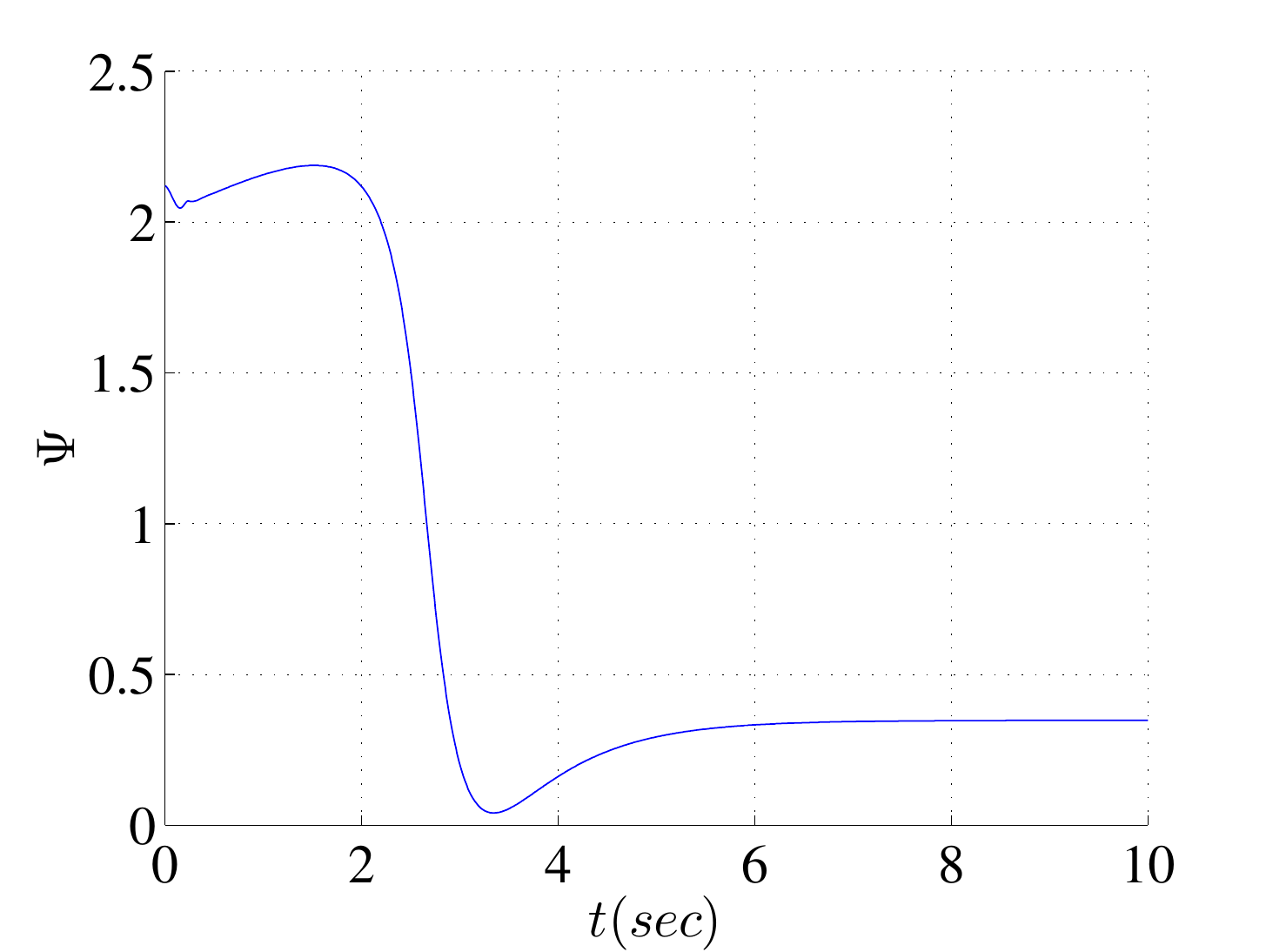} 
		\caption{Configuration error \( \Psi \)} \label{fig:Psi_con} 
	\end{subfigure}
	\caption{Attitude stabilization without adaptive update law}
	\label{fig:con} 
\end{figure}
\begin{figure} 
	\centering 
	\begin{subfigure}[htbp]{0.5\columnwidth} 
		\includegraphics[width=\columnwidth]{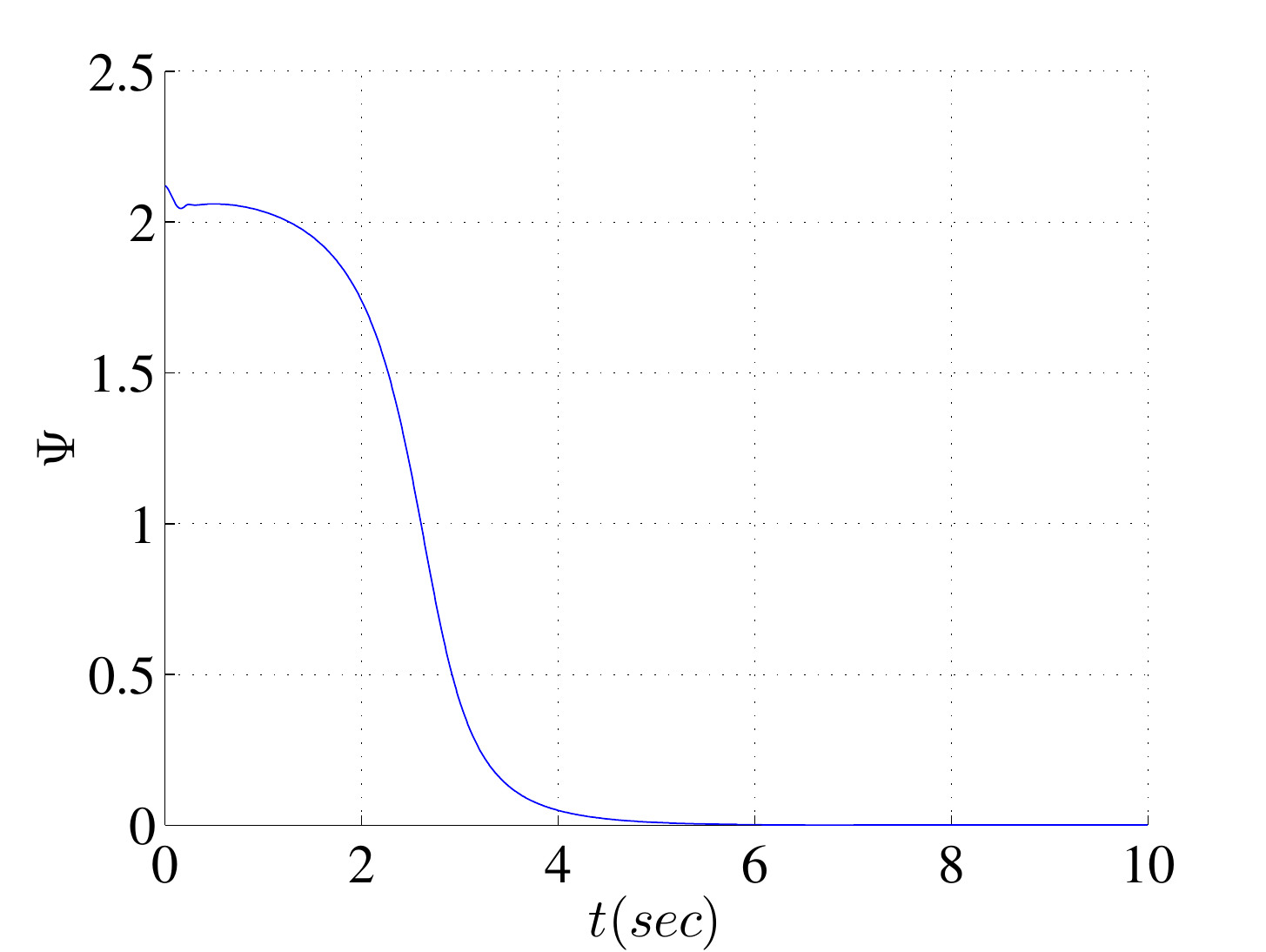} 
		\caption{Configuration error \( \Psi \)} \label{fig:Psi_adapt} 
	\end{subfigure}~ 
	\begin{subfigure}[htbp]{0.5\columnwidth} 
		\includegraphics[width=\columnwidth]{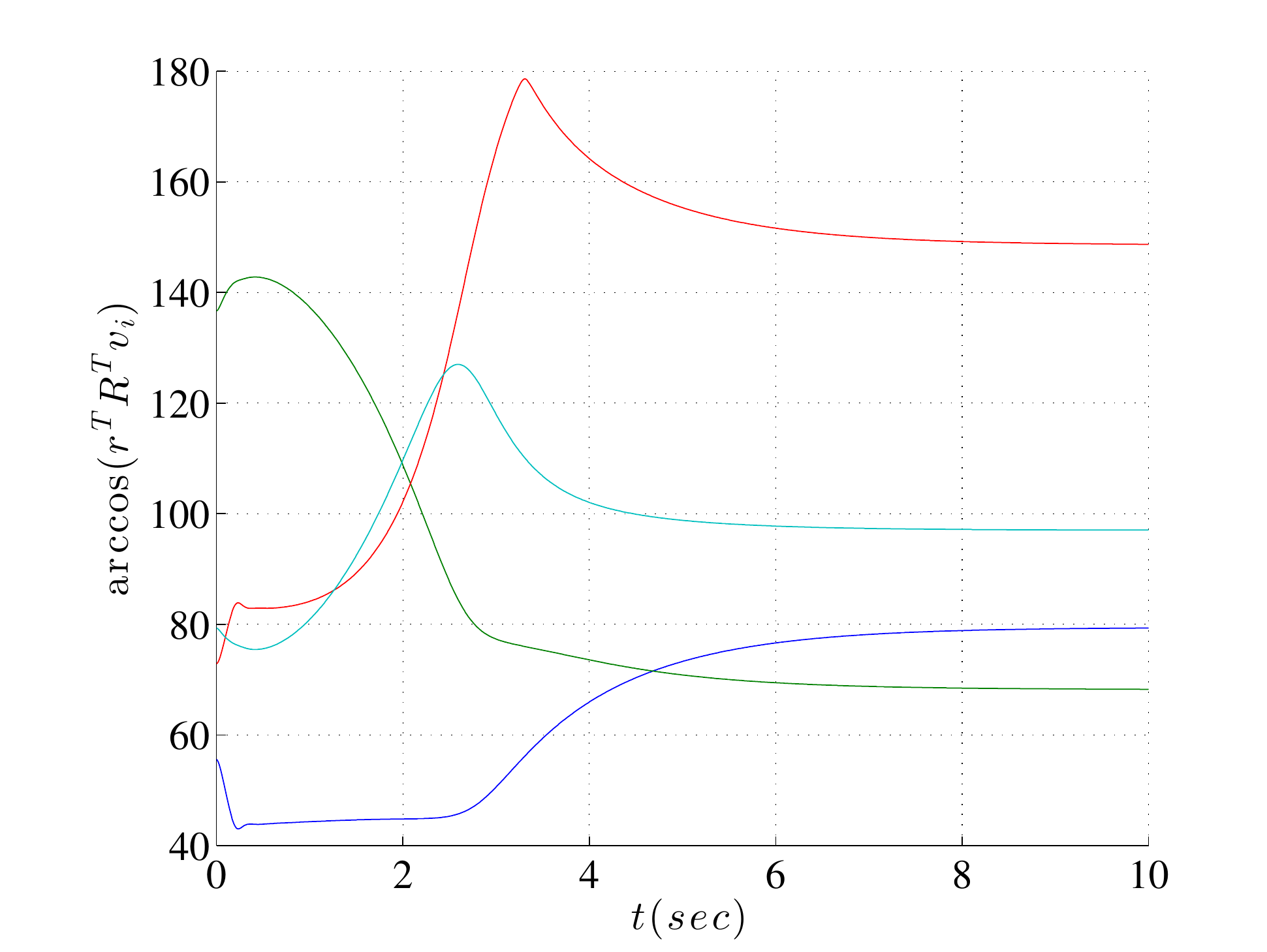} 
		\caption{Angle to constraints} \label{fig:con_angles} 
	\end{subfigure}
	
	\centering
	\begin{subfigure}[htbp]{0.5\columnwidth} 
		\includegraphics[width=\columnwidth]{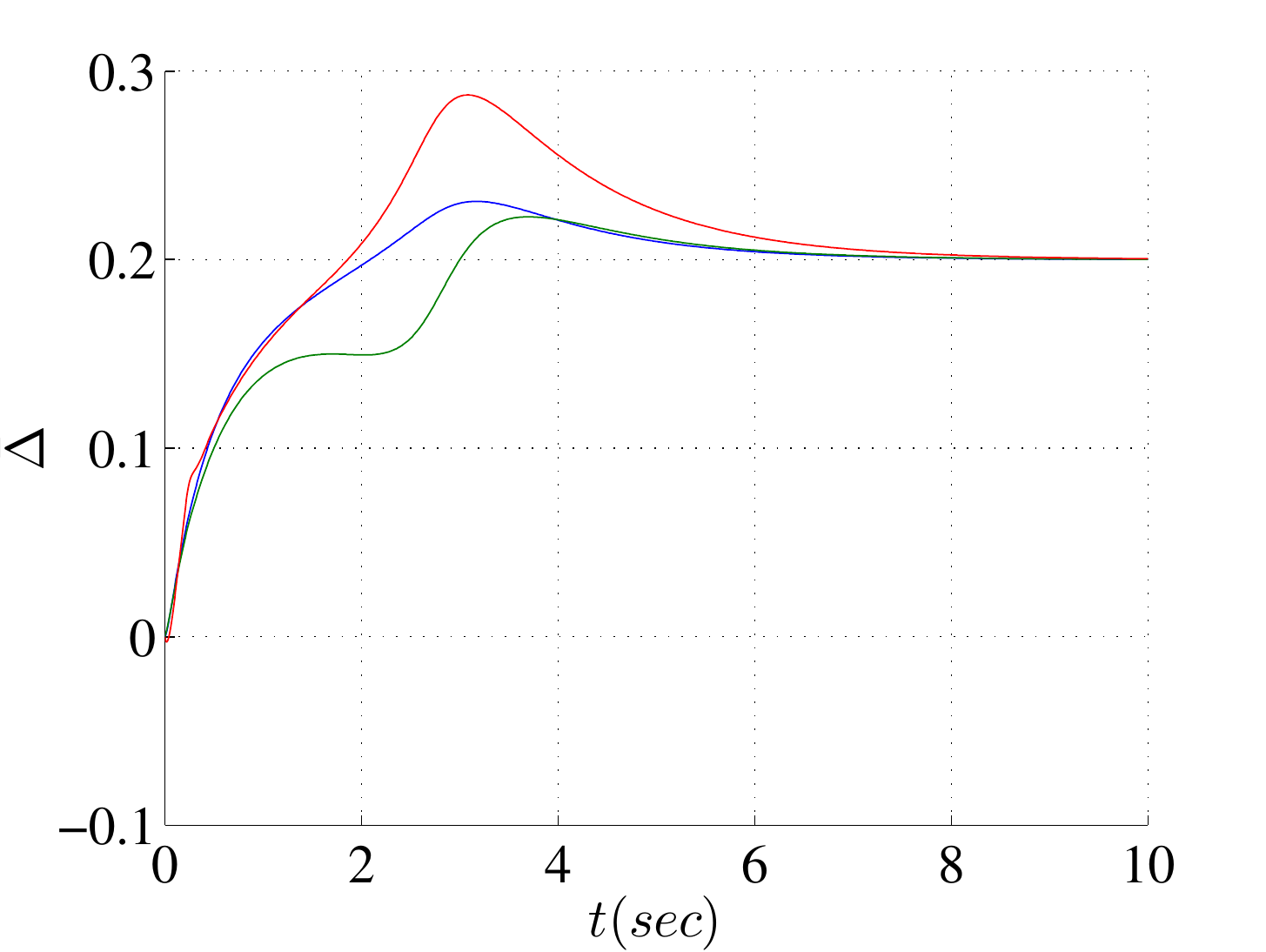} 
		\caption{Disturbance estimate \( \bar \Delta \)} \label{fig:delta_adapt} 
	\end{subfigure}~
	\begin{subfigure}[htbp]{0.5\columnwidth} 
		\includegraphics[trim=9cm 2.5cm 0 0,clip,width=1.29\columnwidth]{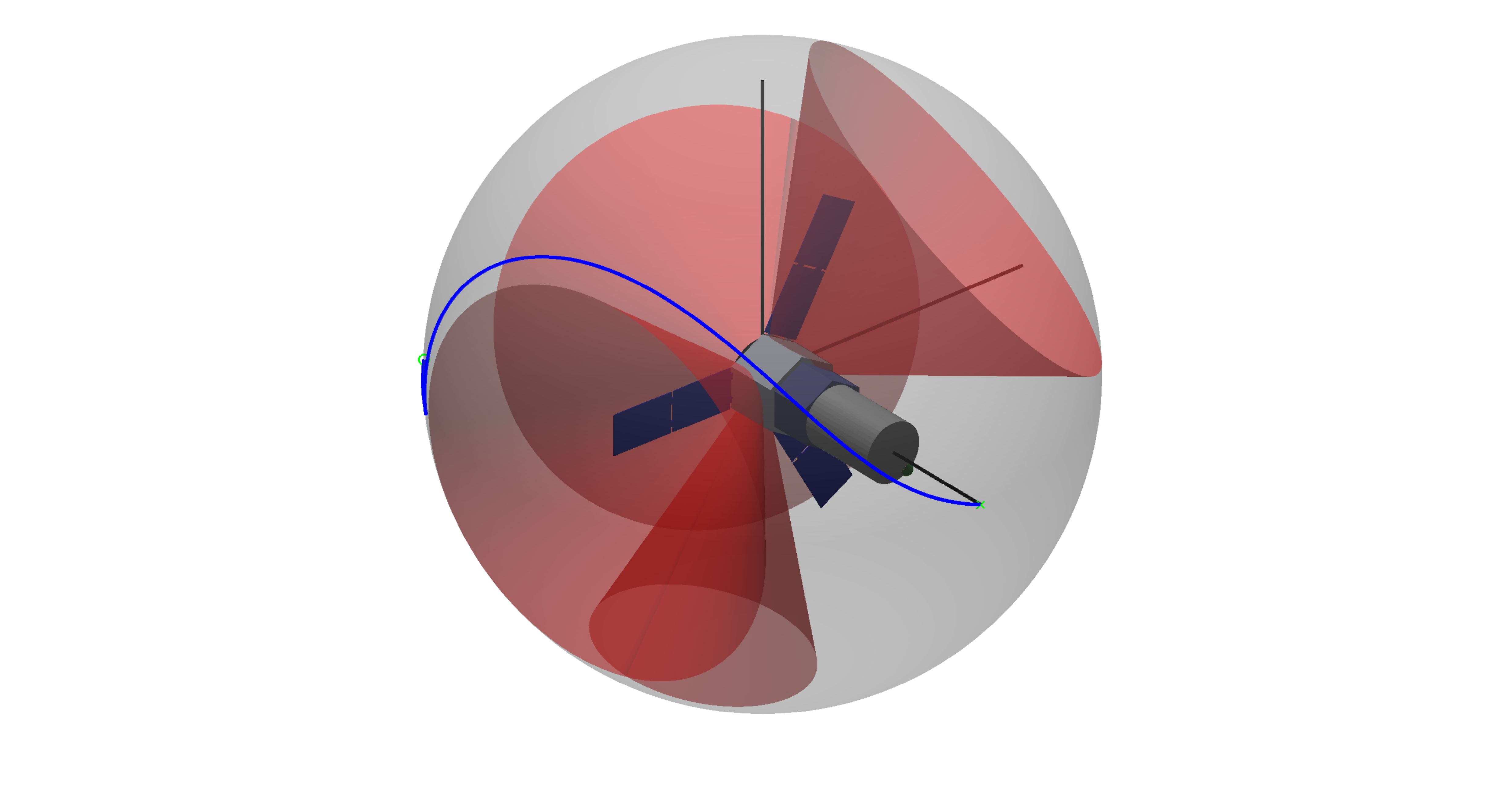} 
		\caption{Attitude trajectory} \label{fig:cad_adapt} 
	\end{subfigure}
	\caption{Attitude stabilization with adaptive update law}
	\label{fig:adapt} 
\end{figure}

Simulation results without the adaptive update law are shown in~\cref{fig:con}.
Without the update law, the system does not achieve zero steady state error. 
\Cref{fig:Psi_con} shows that the configuration error function does not converge to zero and there exist steady state errors.
\Cref{fig:adapt} shows the results with the addition of the adaptive update law.
The addition of the adaptive update law allows the system to converge to the desired attitude in the presence of constraints.
The path of the body fixed sensor in the inertial frame, namely \( R r \), is illustrated in~\cref{fig:cad_adapt}.
The initial attitude is represented with the green circle while the final attitude is marked with a green \(\times\).
The inequality constraints from~\Cref{tab:constraints} are depicted as red cones, where the cone half angle is \( \theta \).
The control system is able to asymptotically converge to zero attitude error.
\Cref{fig:con_angles} shows that the angle \( \arccos(r^T R^T v_i) \) between the body fixed sensor and each constraint is satisfied for the entire maneuver.
In addition, the estimate of the disturbance converges to the the true value as shown in~\cref{fig:delta_adapt}.

Both control system are able to automatically avoid the constrained regions. 
In addition, these results show that it is straightforward to incorporate an arbitrary amount of large constraints.
In spite of this challenging configuration space the proposed control system offers a simple method of avoiding constrained regions.
These closed-loop feedback results are computed in real time and offer a significant advantage over typical open-loop planning methods.
These results show that the proposed geometric adaptive approach is critical to attitude stabilization in the presence of state constraints and disturbances.
\section{Experiment on Hexrotor UAV}
A hexrotor unmanned aerial vehicle (UAV) has been developed at the Flight Dynamics and Controls Laboratory (FDCL) at the George Washington University~\cite{kaufman2014}.
The UAV is composed of three pairs of counter-rotating propellers. 
The propeller pairs of the hexrotor are angled relative to one another to allow for a fully actuated rigid body.

\begin{figure}
	\centering
	\includegraphics[width = 0.7\columnwidth]{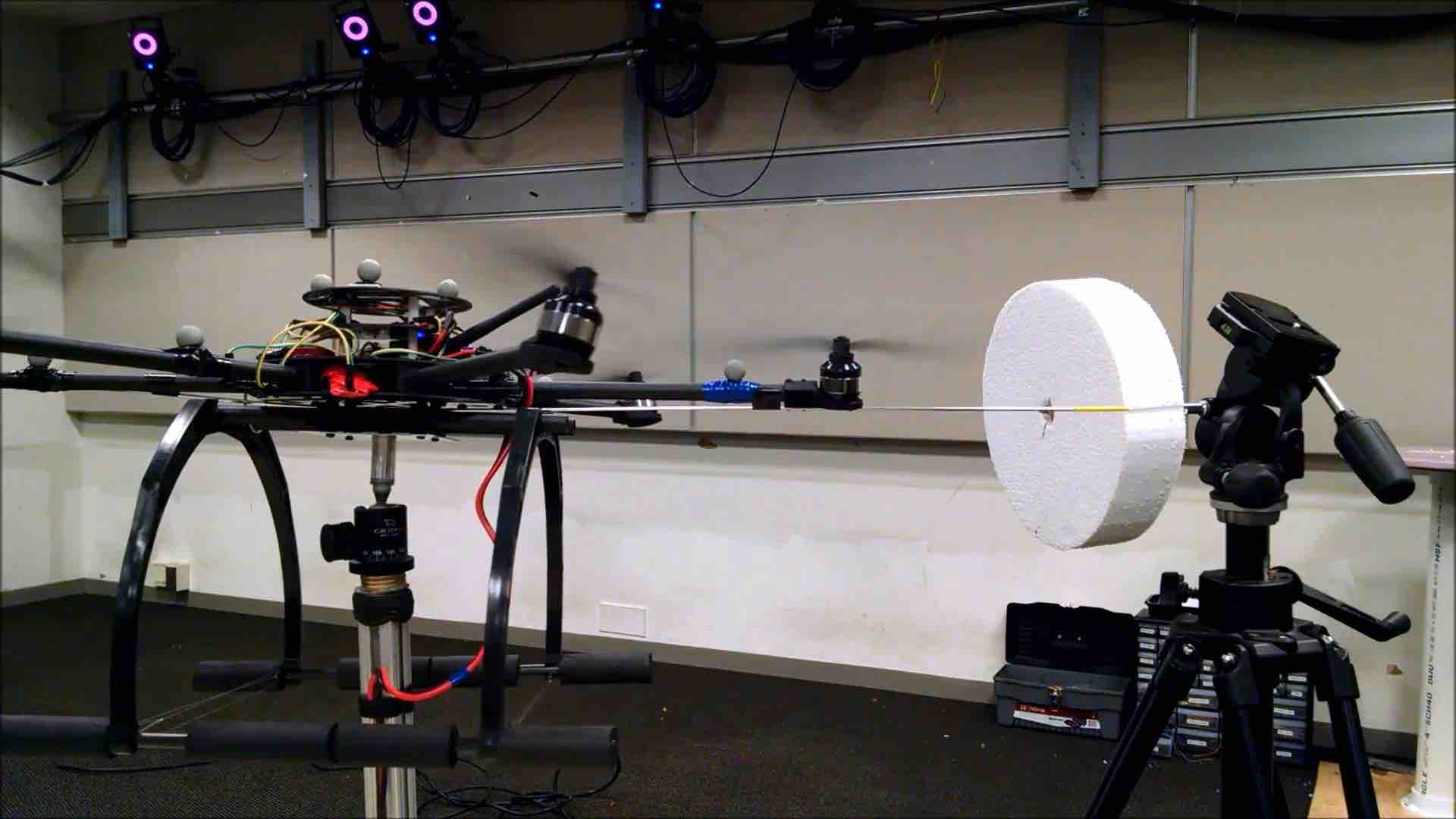}
	\caption{Attitude control testbed~\label{fig:hexrotor}}
\end{figure}

The hexrotor UAV, shown in~\cref{fig:hexrotor}, is composed of the following hardware:
\begin{itemize}
	\item Onboard ODROID XU3 computer module.
	\item VectorNav VN100 IMU operating via TTL serial 
	\item BLDC motors with BL-Ctrl-2.0 ESC via I2C.
	\item Position and attitude over WiFi (TCP) communication from Vicon motion capture system.
	\item Commands sent over WiFi to  onboard controller. 
\end{itemize}
In order to constrain the motion and test only the attitude dynamics we attach the hexrotor to a spherical joint.
The center of rotation is below the center of gravity of the hexrotor.
As a result, there is a destabilizing gravitational moment and the resulting attitude dynamics are similar to an inverted pendulum model.
We augment the control input in~\cref{eqn:adaptive_control} with an additional term to negate the effect of the gravitational moment.

A sensor pointing direction is defined in the body frame of the hexrotor as \( r = [1,0,0]^T \).
We define an obstacle in the inertial frame as \( v = [\frac{1}{\sqrt{2}}, \frac{1}{\sqrt{2}}, 0]^T \) with \( \theta = \ang{12} \).
An initial state is defined as \(R(0) = \exp( \frac{\pi}{2} \hat{e}_3) \), while the desired state is \(R_d =I \).
This results in the UAV performing a \ang{90} yaw rotation about the vertical axis of the spherical joint and the constrained region is on the shortest path connecting $R_0$ and $R_d$. 
The attitude control system is identical to the one presented in~\Cref{prop:adaptive_control} with the exception of a gravity moment term and the following parameters: \(k_R = 0.4, k_\Omega = 0.7 ,c = 0.1 , \alpha = 8 \text{ and } k_\Delta = 0.05\).
\begin{figure} 
	\centering 
	\begin{subfigure}[htbp]{0.5\columnwidth} 
		\includegraphics[width=\columnwidth]{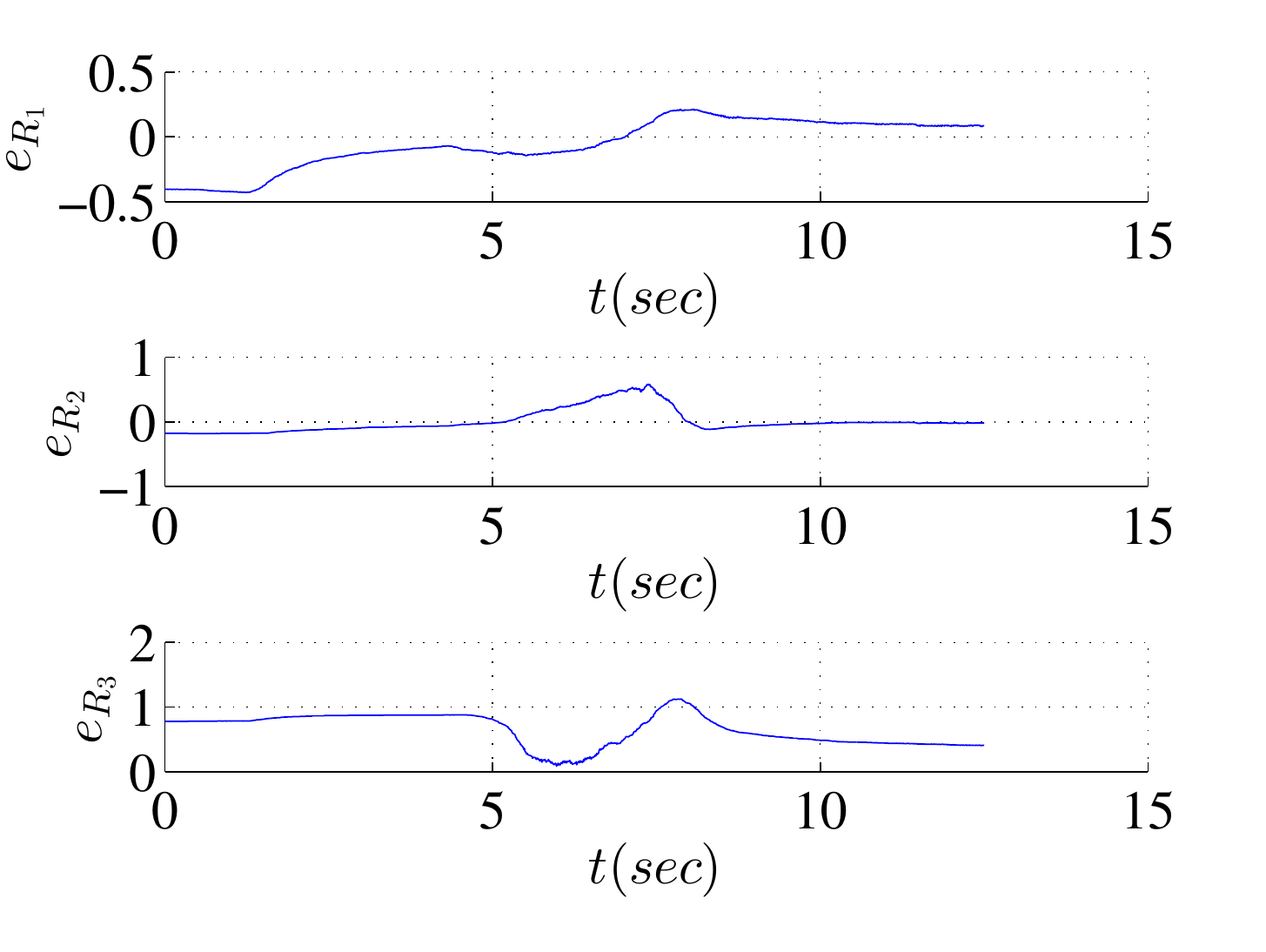} 
		\caption{Attitude error vector \(e_R\) } \label{fig:eR_exp} 
	\end{subfigure}~ 
	\begin{subfigure}[htbp]{0.5\columnwidth} 
		\includegraphics[width=\columnwidth]{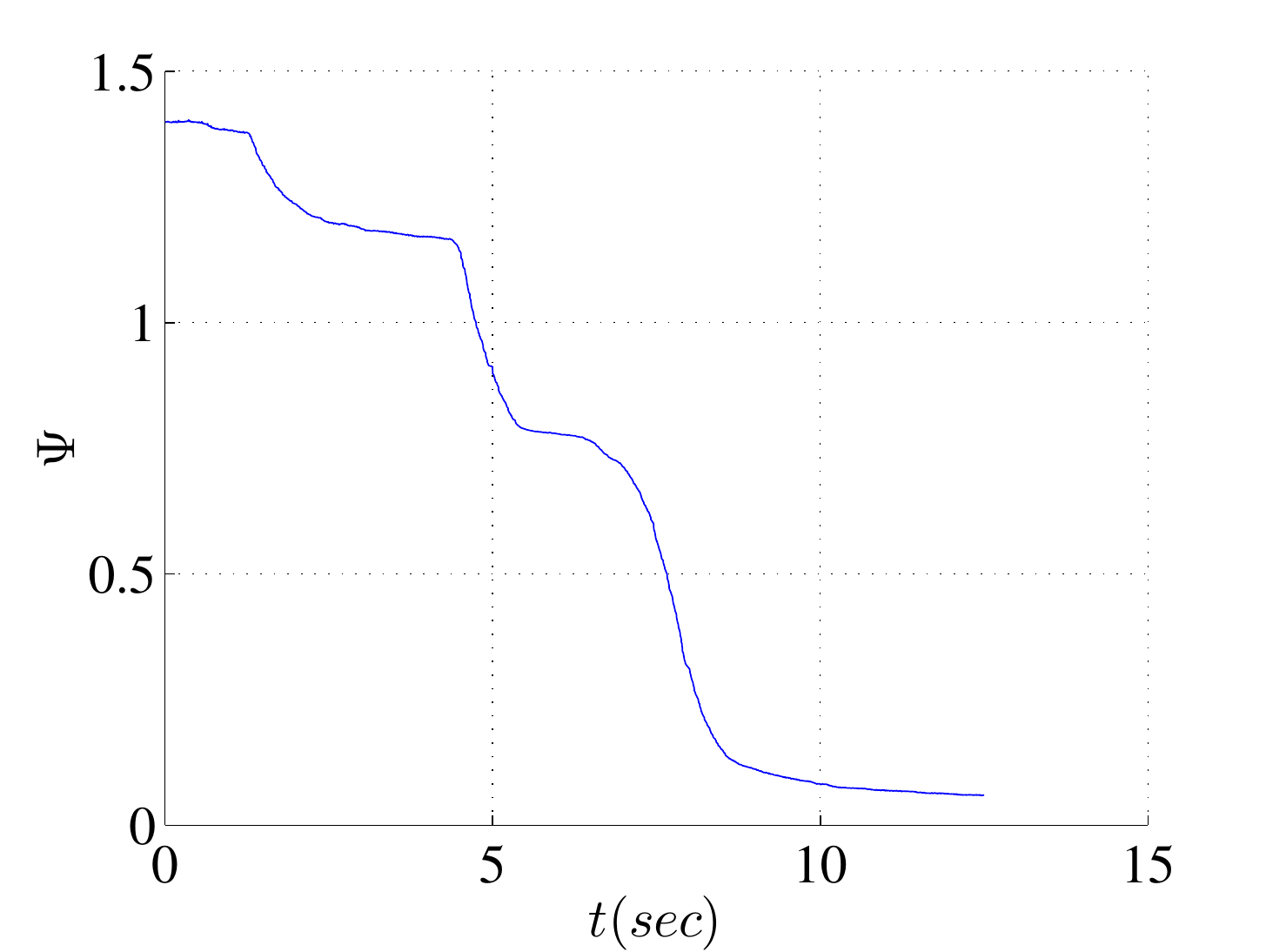} 
		\caption{Configuration error \( \Psi \)} \label{fig:Psi_exp} 
	\end{subfigure}
	
	\begin{subfigure}[htbp]{0.5\columnwidth} 
		\includegraphics[width=\columnwidth]{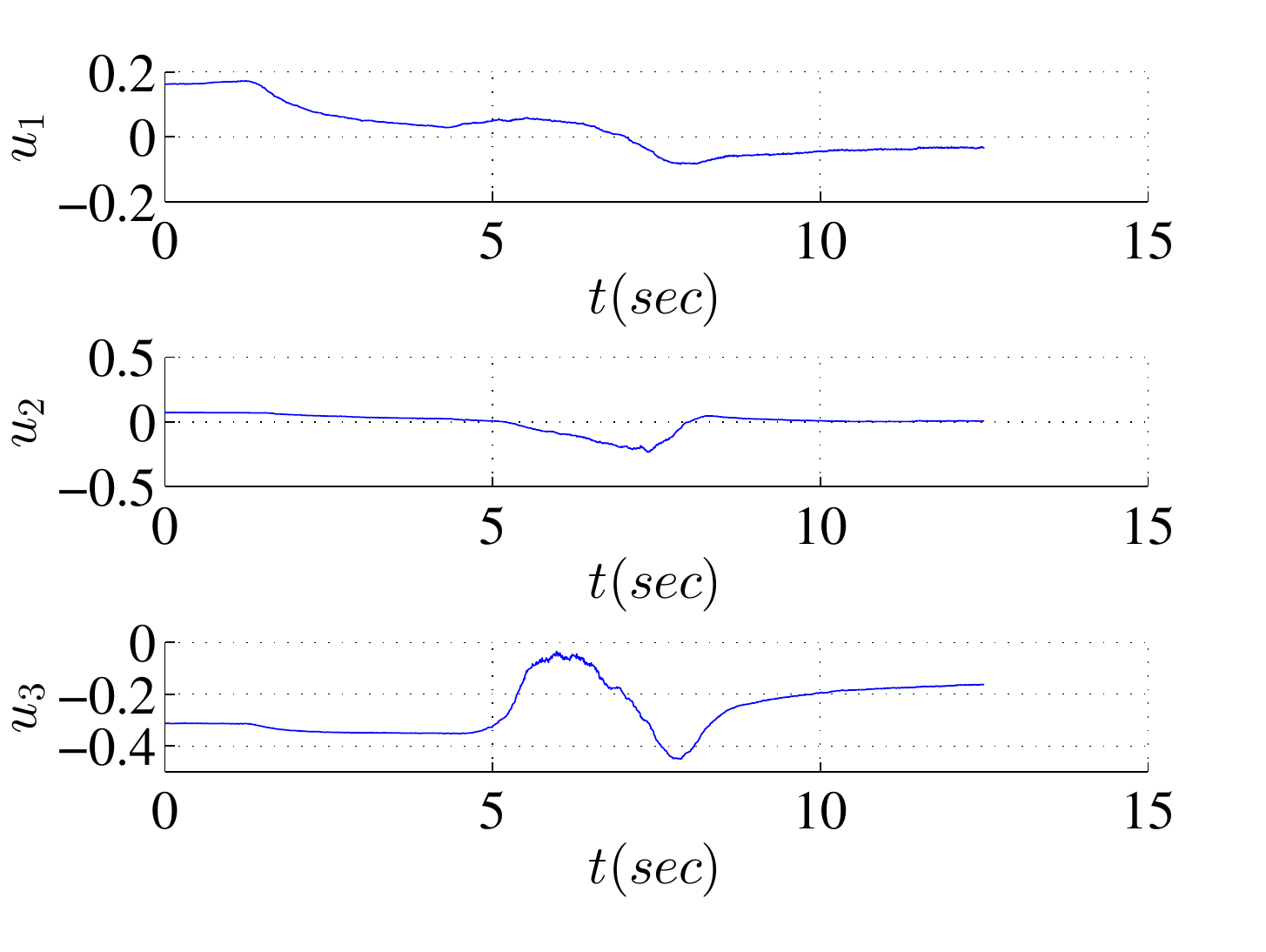} 
		\caption{Control input \( u\)} \label{fig:u_exp} 
	\end{subfigure}~
	\begin{subfigure}[htbp]{0.5\columnwidth} 
		\includegraphics[width=\columnwidth]{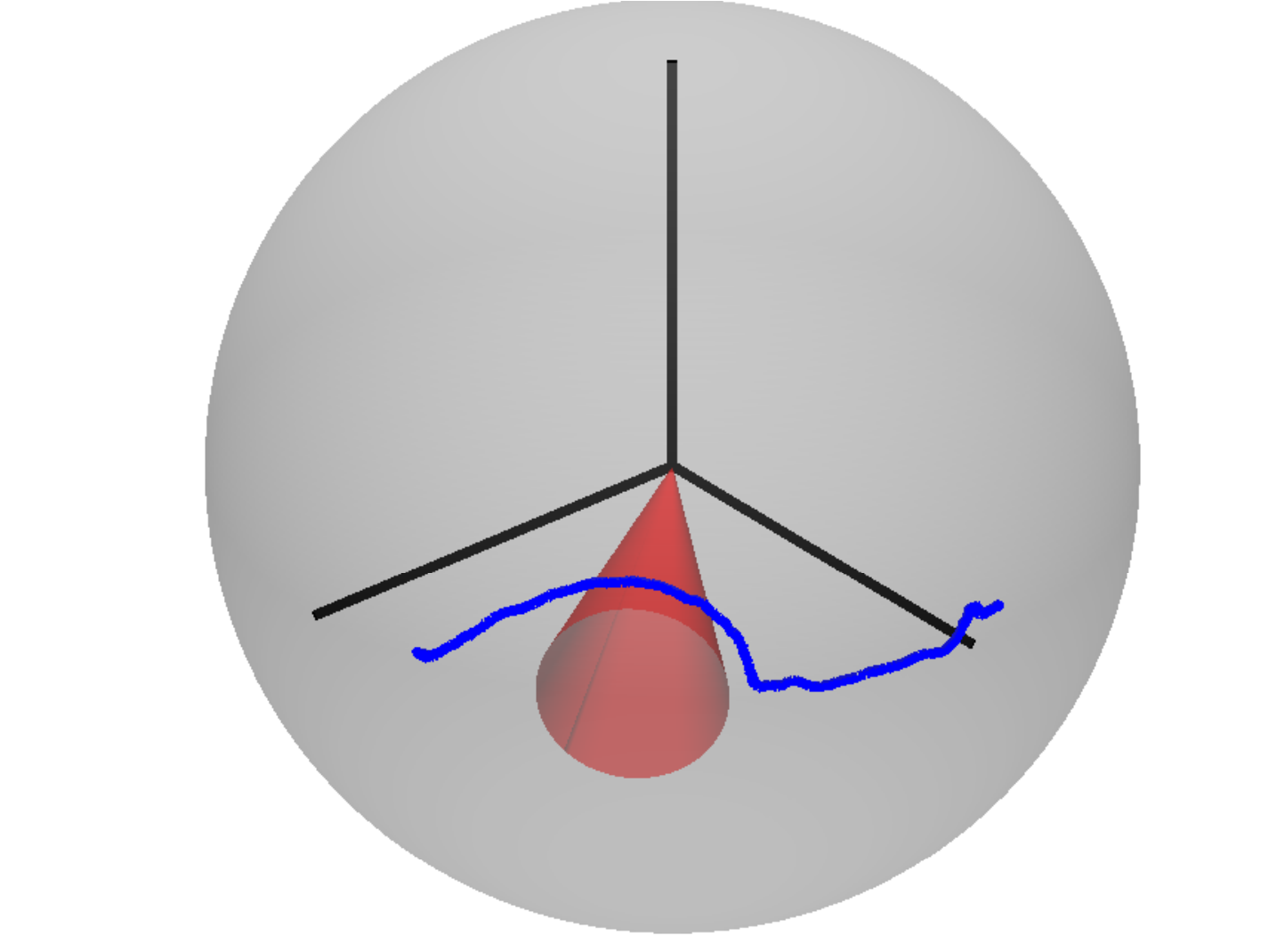} 
		\caption{Attitude Trajectory} \label{fig:traj_exp} 
	\end{subfigure}
	\caption{Constrained Attitude stabilization experiment}
	\label{fig:exp} 
\end{figure}

The experimental results are shown in~\Cref{fig:exp}.
In order to maneuver the system ``close" to the constrained zone we utilize several intermediary set points on either side of the obstacle.
From the initial attitude the hexrotor rotates to the first set point, pauses, and then continues around the obstacle to the second set point before continuing toward the desired attitude.
As a result this creates the stepped behavior of the configuration error history as shown in~\cref{fig:Psi_exp}.

The brushless motors of the hexrotor allow for large control inputs which are critical to enable aggressive maneuvers. 
When constrained to the spherical joint the hexrotor is capable of performing responsive attitude changes with high angular velocities.
In addition, The on-board control and motion capture system operate at a discrete interval of approximately \SI{100}{\hertz}.
It is possible for the system to violate the constraint between these discrete steps and cause numerical exceptions within the embedded software. 
As a result, conservative control gains are chosen to ensure the hexrotor operates in a sedate manner and to allow sufficient time for the measurement and control software to operate.

There exist several sources of error in the experimental setup.
The motion capture system uses a series of optical sensors to determine the relative position of several tracking markers. 
These markers as well as the cameras must remain fixed to ensure accurate attitude measurement.
In addition, the spherical joint is not fixed at the center of mass but is instead offset due to the physical structure of the hexrotor.
As a result a disturbance moment is induced on the resultant motion.

This results in a small steady state error in the vicinity of the desired attitude. 
Over time~\cref{eqn:delta_dot} will remain non-zero while \( e_R \neq 0 \).
This will cause an increase in control input until the steady-state error is reduced. 
Further tuning of the control gains would enable a faster response and a reduced settling time.

The hexrotor avoids the constrained region illustrated by the circular cone in~\cref{fig:traj_exp}, by rotating around the boundary of the constraint. 
This verifies that the proposed control system exhibits the desired performance in the experimental setting as well. 
A video clip showing the attitude maneuver is available \url{https://youtu.be/dsmAbwQram4}.
\section{Conclusions}\label{sec:conclusions}
We have developed a geometric adaptive control system which incorporates state inequality constraints on \(\SO\).
The presented control system is developed directly on \(\SO\) and it avoids singularities and ambiguities that are inherent to attitude parameterizations.
The attitude configuration error is augmented with a barrier function to avoid the constrained region, and an adaptive control law is proposed to cancel the effects of uncertainties. 
We show the stability of the  proposed control system through a rigorous mathematical analysis.
In addition, we have demonstrated the control system via numerical simulation and hardware experiments on a hexrotor UAV.
A novel feature of this control is that it is computed autonomously on-board the UAV.
This is in contrast to many state constrained attitude control systems which require an a priori attitude trajectory to be calculated. 
The presented method is simple, efficient and ideal for hardware implementation on embedded systems.

\appendix
\subsection{Proof of~\Cref{prop:config_error}}\label{proof:config_error}
To prove~\cref{item:prop_psi_psd} we note that~\cref{eqn:A} is a positive definite function about \( R = R_d \)~\cite{bullo2004}.
The constraint angle is assumed \( \ang{0} \leq \theta \leq \ang{90} \) such that \( 0 \leq \cos \theta \).
The term \( r^T R^T v \) represents the cosine of the angle between the body fixed vector \( r \) and the inertial vector \( v \). 
It follows that
\begin{align*}
	0 \leq  \frac{\cos \theta -  r^T R^T v}{1 + \cos \theta} \leq 1 ,
\end{align*}
for all \( R \in \SO \). 
As a result, its negative logarithm is always positive and from~\cref{eqn:B}, \(1 < B\).
The error function \( \Psi = A B \) is composed of two positive terms and is therefore also positive definite, and it is minimized at \( R = R_d \).

Next, we consider~\cref{item:prop_era}.
The variation of~\cref{eqn:A} is taken with respect to \( \delta R = R \hat \eta \) as
\begin{align*}
	\dirDiff{A}{R} &= \eta \cdot \frac{1}{2} \left( G R_d^T R - R^T R_d G\right)^\vee ,
\end{align*}
where we used~\cref{eqn:hat1}.

A straightforward application of the chain and product rules of differentiation allows us to show~\cref{item:prop_erb} as
\begin{align*}
	\dirDiff{B}{R} &=  \eta \cdot \frac{ - \left( R^T v\right)^\vee r}{\alpha \left(\cos \theta - r^T R^T v \right)} ,
\end{align*}
where the scalar triple product~\cref{eqn:STP} was used.

The critical points of \( e_{R_A} \) are derived in~\cite{bullo2004}.
There are four critical points of \( e_{R_A} \), the desired attitude \( R_d \) as well as rotations about each body fixed axis by \ang{180}.
The repulsive error vector \( e_{R_{B}} \) is zero only when the numerator \( \parenth{R^T v}^\wedge r = 0 \). 
This condition only occurs if the desired attitude results in the body fixed vector \( r \) becoming aligned with \(R^T v \) while simultaneously satisfying \(\braces{R_d} \cup \braces{R_d \exp(\pi \hat{s}} \) for \( s \in \braces{e_1, e_2, e_3} \).
Since we assume the system will not operate in violation of the constraints, the addition of the barrier function does not add additional critical points to the control system.
The desired equilibrium is \( e_R = 0 \) and \( A = 0\).
The proof of \( \norm{e_{R_A}} \) given by~\cref{item:prop_era_upbound} is available in~\cite{LeeITCST13}.
\subsection{Proof of~\Cref{prop:error_dyn}}\label{proof:error_dyn}
From the kinematics~\cref{eqn:Rdot} and noting that \( \dot{R}_d = 0 \) the time derivative of \( R_d^T R \) is given as
\begin{gather*}
	\diff{}{t} \parenth{R_d^T R} = R_d^T R \hat{e}_\Omega .
\end{gather*}
Applying this to the time derivative of~\cref{eqn:A} gives
\begin{gather*}
	\diff{}{t} (A) = -\frac{1}{2} \tr{G R_d^T R \hat{e}_\Omega} .
\end{gather*}
Applying~\cref{eqn:hat1} into this shows~\cref{eqn:A_dot}.
Next, the time derivative of the repulsive error function is given by
\begin{gather*}
	\diff{}{t} (B) = \frac{r^T \parenth{\hat{\Omega} R^T} v}{\alpha \parenth{r^T R^T v - \cos \theta}} .
\end{gather*}
Using the scalar triple product, given by~\cref{eqn:STP}, one can reduce this to~\cref{eqn:B_dot}.
The time derivative of the attractive attitude error vector, \( e_{R_A} \), is given by
\begin{gather*}
	\diff{}{t} ( e_{R_A}) = \frac{1}{2} \parenth{\hat{e}_\Omega R^T R_d G + (R^T R_d G)^T \hat{e}_\Omega}^\vee .
\end{gather*}
Using the hat map property given in~\cref{eqn:xAAx} this is further reduced to~\cref{eqn:eRA_dot,eqn:E}.

We take the time derivative of the repulsive attitude error vector, \( e_{R_B} \), as
\begin{gather*}
	\diff{}{t}( e_{R_B} )= a \Omega v^T R r - a R^T v \Omega^T r + b R^T \hat{v} R r ,
\end{gather*}
with \( a \in \R \) and \( b \in \R\) given by 
\begin{gather*}
	a = \bracket{\alpha \parenth{r^T R^T v - \cos \theta}}^{-1} , \,
	b = \frac{r^T \hat{\Omega} R^T v}{\alpha \parenth{r^T R^T v - \cos \theta}^2} .
\end{gather*}
Using the scalar triple product from~\cref{eqn:STP} as \( r \cdot \Omega \times \parenth{R^T v} = \parenth{R^T v} \cdot r \times \Omega \) gives~\cref{eqn:eRB_dot,eqn:F}.

We show the time derivative of the configuration error function as
\begin{gather*}
	\diff{}{t} (\Psi) = \dot{A} B + A \dot{B} .
\end{gather*}
A straightforward substitution of~\cref{eqn:A_dot,eqn:B_dot,eqn:A,eqn:B} into this and appplying~\cref{eqn:eR} shows~\cref{eqn:psi_dot}.
We show~\cref{eqn:eW_dot} by rearranging~\cref{eqn:Wdot} as 
\begin{align*}
	\diff{}{t} e_\Omega = \dot{\Omega} = J^{-1} \parenth{u - \Omega \times J \Omega + W(R,\Omega) \Delta } .
\end{align*}

\subsection{Proof of~\Cref{prop:att_control}}\label{proof:att_control}
Consider the following Lyapunov function:
\begin{gather}
	\mathcal{V} = \frac{1}{2} e_\Omega \cdot J e_\Omega + k_R \Psi(R,R_d) . \label{eqn:v_nodist}
\end{gather}
From~\cref{item:prop_psi_psd} of~\Cref{prop:config_error}, \(\mathcal{V} \geq 0 \).
Using~\cref{eqn:eW_dot,eqn:psi_dot} with \( \Delta = 0 \), the time derivative of \( \mathcal{V} \) is given by
\begin{align}
	\dot{\mathcal{V}} &= -k_\Omega \norm{e_\Omega}^2 . \label{eqn:vdot_nodist}
\end{align}
Since \( \mathcal{V} \) is positive definite and \( \dot{\mathcal{V}} \) is negative semi-definite, the zero equilibrium point \( e_R, e_\Omega \) is stable in the sense of Lyapunov. 
This also implies \( \lim_{t\to\infty} \norm{e_\Omega} = 0 \) and \( \norm{e_R} \) is uniformly bounded, as the Lyapunov function is non-increasing. From \refeqn{eRA_dot} and \refeqn{eRB_dot}, $\lim_{t\to\infty} \dot e_R =0$. 
One can show that \( \norm{\ddot{e}_R} \) is bounded.
From Barbalat's Lemma, it follows \( \lim_{t\to\infty}\norm{\dot{e}_R} = 0 \)~\cite[Lemma 8.2]{khalil1996}. 
Therefore, the equilibrium is asymptotically stable. 
	
Furthermore, since \( \dot{\mathcal{V}} \leq 0 \) the Lyapunov function is uniformly bounded which implies 
\begin{align*}
	\Psi(R(t)) \leq \mathcal{V}(t) \leq \mathcal{V}(0) .
\end{align*}
In addition, the logarithmic term in~\cref{eqn:B} ensures \( \Psi(R) \to \infty \) as \( r^T R^T v \to \cos \theta \).
Therefore, the inequality constraint is always satisfied given that the desired equilibrium lies in the feasible set.
	
\subsection{Proof of~\Cref{prop:eR_dot_bound}}\label{proof:eR_dot_bound}
	The selected domain ensures that the configuration error function is bounded \( \Psi < \psi \).
	This implies that that both \( A(R) \) and \( B(R) \) are bounded by constants \( c_A c_B < \psi < h_1\).
	Furthermore, since \( \norm{B} > 1 \) this ensures that \( c_A, c_B < \psi\) and shows~\cref{eqn:AB_bound}.
	
	Next, we show~~\cref{eqn:E_bound,eqn:F_bound} using the Frobenius norm.
	The Frobenius norm \( \norm{E}_F \) is given in~\cite{LeeITCST13} as
	\begin{gather*}
		\norm{E}_F = \sqrt{\tr{E^T E}} = \frac{1}{2} \sqrt{\tr{G^2} + \tr{R^T R_d G}^2} .
	\end{gather*}
	Applying Rodrigues' formula and the Matlab symbolic toolbox, this is simplified to
	\begin{gather*}
		\norm{E}^2_F \leq \frac{1}{4} \parenth{\tr{G^2} + \tr{G}^2} \leq \frac{1}{2} \tr{G}^2 ,
	\end{gather*}
	which shows~\cref{eqn:E_bound}, since \( \norm{E} \leq \norm{E}_F \).
	
	To show~\cref{eqn:F_bound}, we apply the Frobenius norm \( \norm{F}_F \):
	\begin{gather*}
		\norm{F}_F = \frac{1}{\alpha ^2 \parenth{r^T R^T v - \cos \theta}^2} \left[\tr{a^T a} - 2 \tr{a^T b} \right. \\
		\left.+ 2 \tr{a^T c } + \tr{b^T b}  - 2 \tr{b^T c} + \tr{c^T c}\right] .
	\end{gather*}
	where the terms \( a, b, \text{ and } c \) are given by
	\begin{gather*}
		a = r^T R r I , \quad	b = R^T v r^T , \quad c = \frac{R^T \hat{v} R r v^T R \hat{r}}{r^T R^T v - \cos \theta}.
	\end{gather*}
	A straightforward computation of \( a^T a \) shows that
	\begin{gather*}
		\tr{a^T a} = \parenth{v^T R r}^2 \tr{I} \leq 3 \beta^2 ,
	\end{gather*}
	where we used the fact that \( v^T R r = r^T R^T v < \beta \) from our given domain.
	Similarly, one can show that \( \tr{a^T b} \) is equivalent to
	\begin{gather*}
		\tr{a^T b} = v^T R r \tr{R^T v r^T} = \parenth{v^T R r}^2 \leq \beta^2 ,
	\end{gather*} 
	where we used the fact that \( \tr{x y^T} = x^T y \).
	The product \( \tr{a^T c} \) is given by
	\begin{gather*}
		\tr{a^T c} = \frac{v^T R r}{r^T R^T v - \cos \theta} \tr{\parenth{R^T v}^\vee \parenth{r v^T R} \hat{r} } ,
	\end{gather*}
	where we used the hat map property~\cref{eqn:RxR}.
	One can show that \(\mathrm{tr}[a^T c] \leq 0 \) over the range \( -1 \leq v^T R r \leq \cos \theta \). 
	Next, \( \tr{b^T b}\) is equivalent to
	\begin{gather*}
		\tr{b^T b} = \tr{r v^T R R^T v r^T} = 1 ,
	\end{gather*}
	since \( r,v \in \Sph^2\).
	Finally, \( \tr{c^T c} \) is reduced to
	\begin{gather*}
		\tr{c^T c} = \tr{\hat{r} R^T v r^T \bracket{-I + R^T v v^T R} r v^T R \hat{r}} ,
	\end{gather*}
	where we used the fact that \( \hat{x}^2 = - \norm{x}^2 I + x x^T\).
	Expanding and collecting like terms gives
	\begin{gather*}
		\tr{c^T c } = \frac{1 - 2\parenth{v^T R r}^2 + \parenth{v^T R r}^4}{\parenth{r^T R^T v - \cos \theta}^2} . 
	\end{gather*}
	Using the the given domain \( r^T R^T v \leq \beta \) gives the upper bound~\cref{eqn:F_bound}.
	The bound on \( e_{R_A} \) is given in~\cref{eqn:psi_lower_bound} while \( e_{R_B} \) arises from the definition of the cross product \( \norm{a \times b} = \norm{a} \norm{b} \sin \theta \).
	Finally, we can find the upper bound~\cref{eqn:eR_dot} as
	\begin{gather*}
		\norm{\dot{e}_R} \leq \parenth{\norm{B} \norm{E} + 2 \norm{e_{R_A}} \norm{e_{R_B}} + \norm{A}\norm{F}} \norm{e_\Omega} \, .
	\end{gather*}
	Using~\crefrange{eqn:AB_bound}{eqn:eRB_bound} one can define \( H \) in terms of known values.
\subsection{Proof of~\Cref{prop:adaptive_control}}\label{proof:adaptive_control}
	Consider the Lyapunov function \( \mathcal{V} \) given by
	\begin{gather}
		\mathcal{V} = \frac{1}{2} e_\Omega \cdot J e_\Omega + k_R \Psi + c J e_\Omega \cdot e_R + \frac{1}{2 k_\Delta} e_\Delta \cdot e_\Delta , \label{eqn:v_adapt}
	\end{gather}
	over the domain \( D \) in~\cref{eqn:domain}.
	From~\Cref{prop:eR_dot_bound}, the Lyapunov function is bounded in \( D \) by
	\begin{gather}
		\mathcal{V} \leq z^T W z , \label{eqn:v_upper_bound}
	\end{gather}
	where \( e_\Delta = \Delta - \bar{\Delta} \), \( z = [\|e_R\|,\|e_\Omega\|,\|e_\Delta\|]^T\in\R^3 \) and the matrix \(W \in \R^{3 \times 3}\) is given by
	\begin{gather*}
		W = \begin{bmatrix}
			k_R \psi & \frac{1}{2} c \lambda_M & 0 \\
			\frac{1}{2} c \lambda_M & \frac{1}{2} \lambda_M & 0 \\
			0 & 0 & \frac{1}{2 k_\Delta}
		\end{bmatrix} .
	\end{gather*}
	The time derivative of \( \mathcal{V}\) with the control inputs~\cref{eqn:adaptive_control} is given by
	\begin{align}
		\dot{\mathcal{V}} =& - k_\Omega e_\Omega^T e_\Omega + \parenth{e_\Omega + c e_R}^T W e_\Delta - k_R c e_R^T e_R \nonumber\\
		&- k_\Omega c e_R^T e_\Omega + c J e_\Omega^T \dot{e}_R - \frac{1}{k_\Delta} e_\Delta^T \dot{\bar{\Delta}} , \label{eqn:vdot}
	\end{align}
	where we used \( \dot{e}_\Delta = - \dot{\bar{\Delta}} \).
	The terms linearly dependent on \( e_\Delta\) are combined with~\cref{eqn:delta_dot} to yield
	\begin{align*}
		 e_\Delta^T \parenth{W^T \parenth{e_\Omega + c e_R} - \frac{1}{k_\Delta} \dot{\bar{\Delta}}} = 0 . 
	\end{align*}
	Using~\Cref{prop:eR_dot_bound} an upper bound on \( \dot{\mathcal{V}} \) is written as
	\begin{gather*}
		\dot{\mathcal{V}} \leq -\zeta^T M \zeta ,
	\end{gather*}
	where $\zeta=[\|e_R\|,\|e_\Omega\|]\in\R^2$, and the matrix \( M \in \R^{2 \times 2} \) is 
	\begin{gather}
		M = \begin{bmatrix}
			k_R c & \frac{k_\Omega c}{2} \\
			\frac{k_\Omega c}{2} & k_\Omega - c \lambda_M H
		\end{bmatrix} .
	\end{gather}
	If \( c \) is chosen such that~\cref{eqn:c_bound} is satisfied the matrix \( M \) is positive definite.
	This implies that $\dot{\mathcal{V}}$ is negative semidefinite and $\lim_{t\to\infty} \zeta=0$. 
	As the Lyapunov function is non-increasing $z$ is uniformly bounded. 
%
%

                                  
\bibliography{BibMaster,library}

\begin{thebibliography}{10}
\providecommand{\url}[1]{#1}
\csname url@rmstyle\endcsname
\providecommand{\newblock}{\relax}
\providecommand{\bibinfo}[2]{#2}
\providecommand\BIBentrySTDinterwordspacing{\spaceskip=0pt\relax}
\providecommand\BIBentryALTinterwordstretchfactor{4}
\providecommand\BIBentryALTinterwordspacing{\spaceskip=\fontdimen2\font plus
\BIBentryALTinterwordstretchfactor\fontdimen3\font minus
  \fontdimen4\font\relax}
\providecommand\BIBforeignlanguage[2]{{%
\expandafter\ifx\csname l@#1\endcsname\relax
\typeout{** WARNING: IEEEtran.bst: No hyphenation pattern has been}%
\typeout{** loaded for the language `#1'. Using the pattern for}%
\typeout{** the default language instead.}%
\else
\language=\csname l@#1\endcsname
\fi
#2}}

\bibitem{hughes2004}
P.~Hughes, \emph{{Spacecraft Attitude Dynamics}}.\hskip 1em plus 0.5em minus
  0.4em\relax Dover Publications, 2004.

\bibitem{wertz1978}
J.~R. Wertz, \emph{Spacecraft Attitude Determination and Control}.\hskip 1em
  plus 0.5em minus 0.4em\relax Springer, 1978, vol.~73.

\bibitem{bhat2000}
S.~P. Bhat and D.~S. Bernstein, ``A topological obstruction to continuous
  global stabilization of rotational motion and the unwinding phenomenon,''
  \emph{Systems \& Control Letters}, 2000.

\bibitem{shuster1993}
M.~D. Shuster, ``A survey of attitude representations,'' \emph{Navigation},
  vol.~8, no.~9, 1993.

\bibitem{chaturvedi2011a}
N.~Chaturvedi, A.~K. Sanyal, N.~H. McClamroch, \emph{et~al.}, ``Rigid-body
  attitude control,'' \emph{Control Systems, IEEE}, vol.~31, no.~3, pp. 30--51,
  2011.

\bibitem{bullo2004}
F.~Bullo and A.~D. Lewis, \emph{Geometric Control of Mechanical Systems}, ser.
  Texts in Applied Mathematics.\hskip 1em plus 0.5em minus 0.4em\relax New
  York-Heidelberg-Berlin: Springer Verlag, 2004, vol.~49.

\bibitem{MayTeePaCC11}
C.~Mayhew and A.~Teel, ``Synergistic potential functions for hybrid control of
  rigid-body attitude,'' in \emph{Proceedings of the American Control
  Conference}, 2011, pp. 875--880.

\bibitem{LEEITAC15}
T.~Lee, ``Global exponential attitude tracking controls on {$\SO$},''
  \emph{IEEE Transactions on Automatic Control}, vol.~60, no.~10, pp.
  2837--2842, 2015.

\bibitem{hablani1999}
\BIBentryALTinterwordspacing
H.~B. Hablani, ``Attitude commands avoiding bright objects and maintaining
  communication with ground station,'' \emph{Journal of Guidance, Control, and
  Dynamics}, vol.~22, no.~6, pp. 759--767, 2015/09/19 1999. [Online].
  Available: \url{http://dx.doi.org/10.2514/2.4469}
\BIBentrySTDinterwordspacing

\bibitem{frazzoli2001}
E.~Frazzoli, M.~Dahleh, E.~Feron, and R.~Kornfeld, ``A randomized attitude slew
  planning algorithm for autonomous spacecraft,'' in \emph{AIAA Guidance,
  Navigation, and Control Conference and Exhibit, Montreal, Canada}, 2001.

\bibitem{guiggiani2014}
\BIBentryALTinterwordspacing
A.~Guiggiani, I.~Kolmanovsky, P.~Patrinos, and A.~Bemporad, ``Fixed-point
  constrained model predictive control of spacecraft attitude,''
  \emph{arXiv:1411.0479}, 2014. [Online]. Available:
  \url{http://arxiv.org/abs/1411.0479}
\BIBentrySTDinterwordspacing

\bibitem{kalabic2014}
U.~Kalabic, R.~Gupta, S.~Di~Cairano, A.~Bloch, and I.~Kolmanovsky,
  ``Constrained spacecraft attitude control on ${\SO}$ using fast nonlinear
  model predictive control using reference governors and nonlinear model
  predictive control,'' in \emph{American Control Conference (ACC), 2014}, June
  2014, pp. 5586--5593.

\bibitem{gupta2015}
R.~Gupta, U.~Kalabic, S.~Di~Cairano, A.~Bloch, and I.~Kolmanovsky,
  ``Constrained spacecraft attitude control on ${\SO}$ using fast nonlinear
  model predictive control,'' in \emph{American Control Conference (ACC),
  2015}, July 2015, pp. 2980--2986.

\bibitem{rimon1992}
E.~Rimon and D.~E. Koditschek, ``Exact robot navigation using artificial
  potential functions,'' \emph{Robotics and Automation, IEEE Transactions on},
  vol.~8, no.~5, pp. 501--518, 1992.

\bibitem{lee2011b}
U.~Lee and M.~Mesbahi, ``{Spacecraft Reorientation in Presence of Attitude
  Constraints via Logarithmic Barrier Potentials},'' in \emph{Proceedings of
  the American Control Conference}, 2011, pp. {450--455}.

\bibitem{mcinnes1994}
\BIBentryALTinterwordspacing
C.~R. McInnes, ``Large angle slew maneuvers with autonomous sun vector
  avoidance,'' \emph{Journal of Guidance, Control, and Dynamics}, vol.~17,
  no.~4, pp. 875--877, 2015/07/10 1994. [Online]. Available:
  \url{http://dx.doi.org/10.2514/3.21283}
\BIBentrySTDinterwordspacing

\bibitem{LeeITCST13}
T.~Lee, ``Robust adaptive tracking on {$\SO$} with an application to the
  attitude dynamics of a quadrotor {UAV},'' \emph{IEEE Transactions on Control
  Systems Technology}, vol.~21, no.~5, pp. 1924--1930, September 2013.

\bibitem{kaufman2014}
E.~Kaufman, K.~Caldwell, D.~Lee, and T.~Lee, ``Design and development of a
  free-floating hexrotor {UAV} for 6-dof maneuvers,'' in \emph{Proceedings of
  the IEEE Aerospace Conference}, 2014.

\bibitem{khalil1996}
H.~K. Khalil, \emph{Nonlinear Systems}, 3rd~ed.\hskip 1em plus 0.5em minus
  0.4em\relax Prentice Hall New Jersey, 2002.

\end{thebibliography}
\bibliographystyle{IEEEtran}

\end{document}